\documentclass[10pt]{amsart}
\usepackage[T1]{fontenc}
\usepackage{geometry}
\usepackage[latin1]{inputenc}
\usepackage{amsmath}
\usepackage{amsfonts, amssymb, textcomp}
\usepackage[colorlinks=flase, linkcolor=red,urlcolor=green, citecolor=blue]{hyperref}
\usepackage{subeqnarray}

\usepackage{latexsym}
\usepackage{fancyhdr}
\usepackage{longtable}
\usepackage{amsmath, amssymb}
\usepackage{graphicx}

\numberwithin{equation}{section}

\newtheorem{lemma}[subsection]{Lemma}
\newtheorem{theorem}[subsection]{Theorem}
\newtheorem{proposition}[subsection]{Proposition}

\newtheorem{remark}[subsection]{Remark}

\newcommand{\RR}{\mathbb{R}}
\newcommand{\CC}{\mathbb{C}}
\newcommand{\NN}{\mathbb{N}}

\newcommand{\supp}{\operatorname{supp}}

\let\on=\operatorname

\newenvironment{demo}[1]{\par\smallskip\noindent{\bf #1.}}{\par\smallskip}

\title[Surjectivity of the Borel mapping in the mixed setting]
{The surjectivity of the Borel mapping in the mixed setting for ultradifferentiable ramification spaces}

\author[J.~Jim\'{e}nez-Garrido, J.~Sanz, and G.~Schindl]{Javier Jim\'{e}nez-Garrido, Javier Sanz and Gerhard Schindl}

\begin{document}
\begin{abstract}
We consider $r$-ramification ultradifferentiable classes, introduced by J. Schmets and M. Valdivia in order to study the surjectivity of the Borel map, and later on also exploited by the authors in the ultraholomorphic context. We characterize quasianalyticity in such classes, extend the results of Schmets and Valdivia about the image of the Borel map in a mixed ultradifferentiable setting, and obtain a version of the Whitney extension theorem in this framework.

\end{abstract}

\date{\today}

\maketitle
\section{Introduction}


Spaces of ultradifferentiable functions are subclasses of smooth functions on an open set $U\subseteq\RR^d$ having a prescribed growth control on the functions and all their derivatives. Classically, in the literature this growth is measured either by a weight sequence $M$ (e.g. see \cite{Komatsu73}) or a weight function $\omega$ (e.g. see \cite{BraunMeiseTaylor90}) and it is shown that in general both methods yield different classes, see \cite{BonetMeiseMelikhov07}. In both settings one can distinguish between the Roumieu type and Beurling type spaces. In this paper we will exclusively consider classes of both types defined by a weight sequence, respectively denoted by $\mathcal{E}_{\{M\}}(U,\CC)$
and $\mathcal{E}_{(M)}(U,\CC)$ (see Subsection~\ref{section21} for the precise definition), or by $\mathcal{E}_{[M]}(U,\CC)$ when both are referred to at the same time.

Analogously, motivated by solving difference and differential equations, there do also exist classes of ultraholomorphic functions defined in terms of a sequence $M$ (mostly of Roumieu type). The functions are defined on unbounded sectors of the Riemann surface of the logarithm, and in this case the weights $M$ control the growth of the complex derivatives. Closely related are classes of functions admitting asymptotic expansion (again on unbounded sectors of the Riemann surface of the logarithm). For more details and the historic development we refer to the introduction of \cite{injsurj} and the references therein.\vspace{6pt}

An important question in both the ultradifferentiable and ultraholomorphic situation is to establish sufficient and necessary conditions on $M$ under which the Borel map $\mathcal{B}_0$, which assigns to $f$ the infinite jet $(f^{(j)}(0))_{j\in\NN}$, is onto the corresponding sequence spaces $\Lambda_{\{M\}}$ or $\Lambda_{(M)}$ (defined in Subsection~\ref{section33}), see \cite{petzsche}, \cite{Thilliezdivision} and \cite{Sanzflatultraholomorphic}. In the ultradifferentiable setting the so-called {\itshape strong non-quasianalyticity condition} \hyperlink{gamma1}{$(\gamma_1)$} is characterizing this behavior for both types as shown in \cite{petzsche}.

In order to study the surjectivity of the (asymptotic) Borel map (or even to show the existence of continuous linear extension operators, i.e. right inverses of the Borel map) for ultraholomorphic classes defined by $M$ (see~\cite[Section 4]{injsurj} and~\cite[Section 3.3]{dissertationjimenez}, concentrating on the Roumieu type), different (auxiliary) spaces of smooth functions have been introduced and used, whose elements $f$ are having ultradifferentiable growth conditions not for all derivatives $f^{(j)}$, $j\in\NN$, but only for all $f^{(rj)}$, $j\in\NN$, where $r\in\NN_{\ge 1}$ is a (ramification) parameter. We refer to Sections \ref{section33} and \ref{testfunctionspaces} for recalling the definitions of the mentioned spaces in the present work.


The property \hyperlink{gammar}{$(\gamma_r)$}, crucially appearing in our results in~\cite{injsurj}, was introduced by Schmets and Valdivia. It was shown to characterize the surjectivity of the Borel map, and even the existence of an extension map, in Beurling ultradifferentiable $r$-ramified classes~\cite[Proposition 4.3 and Theorem 4.4]{Schmetsvaldivia00}, and to be necessary for the surjectivity of the Borel map in the Roumieu case~\cite[Proposition 5.2]{Schmetsvaldivia00} (in~\cite[Theorem 5.4]{Schmetsvaldivia00} they also show that the existence of an extension map in this case amounts to \hyperlink{gammar}{$(\gamma_r)$} and the restrictive condition $(\beta_2)$ from \cite{petzsche}).

Closely related, one may consider $\mathcal{D}_{[M]}([-1,1])$, the subspace of $\mathcal{E}_{[M]}(\RR,\CC)$ whose elements' support is contained in $[-1,1]$. In \cite{surjectivity} a complete characterization of the fact that $\Lambda_{[M]}\subseteq\mathcal{B}_0(\mathcal{D}_{[N]}([-1,1]))$ was given in a mixed setting between two classes defined by generally different sequences $M$ and $N$, both not having \hyperlink{gamma1}{$(\gamma_1)$} necessarily, see also in \cite{BonetMeiseTaylorSurjectivity} for the weight function setting and in \cite{ChaumatChollet94} working with the more general Whitney jet mapping on compact sets (but assuming more restrictive standard conditions on the weights).\vspace{6pt}

The main aim of this article is to transfer the mixed-setting results from \cite{surjectivity} to these (non-standard) $r$-ramified classes, and the motivation of this question was arising when inspecting the proof of the surjectivity of the asymptotic Borel map in ultraholomorphic classes \cite[Thm. 4.14 $(i)$]{injsurj} for sequences with the property of {\itshape derivation closedness}. Indeed, without this additional assumption on $M$, one obtains a result similar to those in \cite{surjectivity} that we have just described, providing information on the image of the Borel map on a $r$-ramification space in the mixed situation between $M=(M_p)_p$ and the forward-shifted sequence $N=(N_p)_p$ with $N_p=M_{p+1}$.

We expect that a characterization of such situation in terms of some precise growth condition involving $M$, $N$ and $r$, as contained in this paper, will be helpful to obtain mixed results for the (asymptotic) Borel map in ultraholomorphic classes defined by weight sequences, i.e. to transfer the results from \cite{injsurj} to a mixed framework with a control on the loss of regularity.

Related to this aim is the following: In \cite{sectorialextensions} and \cite{sectorialextensions1} the authors have introduced ultraholomorphic classes defined in terms of weight functions $\omega$ and shown partial extension results in this setting. A joint future research will be to completely transfer the results from \cite{injsurj} to the weight function setting, and obtain also in this context some mixed extension procedures.\vspace{6pt}

The paper is organized as follows: First, in Section \ref{section2}, we collect and summarize all necessary notation and conditions on weight sequences $M$, which will be important later on. Moreover we introduce the classical spaces of ultradifferentiable functions defined by weight sequences, the most important ultradifferentiable $r$-ramification function spaces and the corresponding sequence spaces.

In Section \ref{section4} we prove the main result for the Roumieu case (see Theorem \ref{Roumieu-Surjectivitytheorem}), in Section \ref{section5} for the Beurling case (see Theorem \ref{Beurling-Surjectivitytheorem}) which is reduced to the Roumieu case by using a technical result from \cite{ChaumatChollet94}. In the proofs we are following the ideas from \cite{surjectivity} and make necessary changes to deal with the parameter $r$. Even in the case $r=1$, which yields the main statement \cite[Theorem 1.1]{surjectivity}, we are dealing with a slightly more general approach than in \cite{surjectivity} since our assumptions on $M$ and $N$ are weaker, see Remarks \ref{Roumieu-Surjectivitytheoremremark} and \ref{Beurling-Surjectivitytheoremremark}.\vspace{6pt}

In Section \ref{section6} we introduce all further ultradifferentiable $r$-ramification function classes from \cite{Schmetsvaldivia00} and prove that the main results from the previous sections also hold true (see Theorem \ref{finaltheorem}). The special case $M=N$, by assuming some mild standard growth conditions on $M$, shows that property \hyperlink{gammar}{$(\gamma_r)$} is characterizing the surjectivity of the Borel map in all $r$-ramification test function spaces (of both types, so improving the results of Schmets and Valdivia specially in the Roumieu case).


In Section \ref{section62} we show that the new introduced mixed conditions \hyperlink{SV}{$(M,N)_{SV_r}$} and \hyperlink{gammarmix}{$(M,N)_{\gamma_r}$} can also be used to characterize the surjectivity of the more general jet mapping in the mixed weight sequence setting (of Roumieu type) and hence give a Whitney extension theorem involving a ramification parameter $r$.

Finally, in Section \ref{nonquasi}, we prove a full characterization of the (non)quasianalyticity of all ultradifferentiable $r$-ramification function classes by using the classical Denjoy-Carleman theorem for ultradifferentiable classes. Here and in several other questions under consideration in this article the so-called $r$-{\itshape interpolating sequence} introduced in \cite{Schmetsvaldivia00} (see \ref{rinterpolatingsequence}) will play an important role: It helps to reduce the $r$-ramified ultradifferentiable framework to the classical one.

\subsection{General notation}
Throughout this paper we will use the following notation: We denote by $\mathcal{E}$ the class of (complex-valued) smooth functions, $\mathcal{C}^{\omega}$ is the class of all real analytic functions. We will write $\NN_{>0}=\{1,2,\dots\}$ and $\NN=\NN_{>0}\cup\{0\}$. Moreover we put $\RR_{>0}:=\{x\in\RR: x>0\}$, i.e. the set of all positive real numbers. For $k=(k_1,\dots,k_d)\in\NN^d$ we have $|k|=k_1+\dots+k_d$ and $f^{(k)}$ shall denote taking partial derivatives of $f$ with respect to $k=(k_1,\dots,k_d)$. On the class $\mathcal{E}$, for all $U\subseteq\RR^d$ non-empty open, $f\in\mathcal{E}(U,\CC)$ and $x_0\in U$ we introduce the {\itshape Borel map} defined by $\mathcal{B}_{x_0}(f):=(f^{(k)}(x_0))_{k\in\NN^d}$. For given $r\in\NN_{\ge 2}$ we put $\mathcal{B}^r_{x_0}(f):=(f^{(rj)}(x_0))_{j\in\NN^d}$ and for convenience, if $x_0=0$, then we will write $\mathcal{B}$ and $\mathcal{B}^r$ instead of $\mathcal{B}_0$ and $\mathcal{B}^r_0$.\vspace{6pt}

{\itshape Convention:} For convenience we will write $\mathcal{E}_{[M]}$ if either $\mathcal{E}_{\{M\}}$ or $\mathcal{E}_{(M)}$ is considered, but not mixing the cases if statements involve more than one $\mathcal{E}_{[M]}$ symbol (and similarly for all further classes of ultradifferentiable $r$-ramification functions and the associated sequence spaces $\Lambda_{\{M\}}$, $\Lambda_{(M)}$).

\section{Notation and conditions}\label{section2}

\subsection{Weight sequences and classes of ultradifferentiable functions}\label{section21}

$M=(M_p)_p\in\RR_{>0}^{\NN}$ is called a {\itshape weight sequence}, we introduce also $m=(m_p)_p$ defined by $m_p:=\frac{M_p}{p!}$ and $\mu_p:=\frac{M_p}{M_{p-1}}$, $\mu_0:=1$. $M$ is called {\itshape normalized} if $1=M_0\le M_1$ holds true which can always be assumed without loss of generality.\vspace{6pt}

For any weight sequence $M$ and $r>0$ we put $M^{1/r}:=(M_p^{1/r})_{p\in\NN}$.\vspace{6pt}

$(1)$ $M$ is called {\itshape log-convex} if
$$\hypertarget{lc}{(\text{lc})}:\Leftrightarrow\;\forall\;j\in\NN_{>0}:\;M_j^2\le M_{j-1}\cdot M_{j+1},$$
equivalently if $(\mu_p)_p$ is nondecreasing. $M$ is called {\itshape strongly log-convex} if \hyperlink{lc}{$(\text{lc})$} holds for the sequence $m$. If $M$ is log-convex and normalized, then both $M$ and the mapping $j\mapsto(M_j)^{1/j}$ are nondecreasing, e.g. see \cite[Lemma 2.0.4]{diploma}. In this case we get $M_k\ge 1$ for all $k\ge 0$ and
\begin{equation}\label{mucompare}
\forall\;k\in\NN_{>0}:\;\;\;(M_k)^{1/k}\le\mu_k.
\end{equation}

Some properties for weight sequences are very basic and so we introduce for convenience the following set:
$$\hypertarget{LCset}{\mathcal{LC}}:=\{M\in\RR_{>0}^{\NN}:\;M\;\text{is normalized, log-convex},\;\lim_{k\rightarrow\infty}(M_k)^{1/k}=\infty\}.$$

$(2)$ $M$ has {\itshape moderate growth} if
$$\hypertarget{mg}{(\text{mg})}:\Leftrightarrow\exists\;C\ge 1\;\forall\;j,k\in\NN:\;M_{j+k}\le C^{j+k} M_j M_k.$$
A weaker condition is {\itshape derivation closedness},
$$\hypertarget{dc}{(\text{dc})}:\Leftrightarrow\exists\;C\ge 1\;\forall\;j\in\NN:\;M_{j+1}\le C^{j+1} M_j.$$
Note that we can replace in both conditions $M$ by $m$ or by any $M^{1/r}$ (by changing the constant $C$).

$(3)$ $M$ is called {\itshape nonquasianalytic,} if
$$\hypertarget{mnq}{(\text{nq})}:\Leftrightarrow\;\sum_{p=1}^{\infty}\frac{1}{\mu_p}<\infty.$$
If $M$ is log-convex then using {\itshape Carleman's inequality} one can show (for a proof see e.g. \cite[Proposition 4.1.7]{diploma}) that $\sum_{p=1}^{\infty}\frac{M_{p-1}}{M_p}<\infty\Leftrightarrow\sum_{p\ge 1}\frac{1}{(M_p)^{1/p}}<\infty$.

More generally, for arbitrary $r>0$ we call $M$ to be $r$-nonquasianalytic, if
$$\hypertarget{mnqr}{(\text{nq}_r)}:\Leftrightarrow\;\sum_{p\ge 1}\left(\frac{1}{\mu_p}\right)^{1/r}<\infty,$$
so \hyperlink{mnq}{$(\text{nq})$} is precisely \hyperlink{mnqr}{$(\text{nq}_1)$}. Provided that $\mu_p\ge 1$, i.e. $M$ is nondecreasing, we have that \hyperlink{mnqr}{$(\text{nq}_r)$} does imply \hyperlink{mnqr}{$(\text{nq}_s)$} for every $s<r$.

If $M\in\hyperlink{LCset}{\mathcal{LC}}$, then by recalling \cite[Prop. 2.13, Def. 3.3, Thm. 3.4]{Sanzflatultraholomorphic} (see also \cite[p. 145]{injsurj}) we have that the so-called {\itshape exponent of convergence} $\lambda_{(\mu_p)_p}:=\inf\left\{\alpha>0: \sum_{p\ge 1}\left(\frac{1}{\mu_p}\right)^{\alpha}<\infty\right\}$ does coincide with $\limsup_{p\rightarrow\infty}\frac{\log(p)}{\log(\mu_p)}$ and is related to the {\itshape index of quasianalyticity} of $M$, denoted by $\omega(M):=\liminf_{p\rightarrow\infty}\frac{\log(\mu_p)}{\log(p)}$ (also called lower order of $M$), by $\omega(M)=\frac{1}{\lambda_{(\mu_p)_p}}$. Hence in our notation we have now
$$\omega(M)=\sup\{r>0: \hyperlink{mnqr}{(\text{nq}_r)}\}.$$

$(4)$ $M$ has $(\beta_1)$ if
$$\hypertarget{beta1}{(\beta_1)}:\Leftrightarrow\;\exists\;Q\in\NN_{>0}:\;\liminf_{p\rightarrow\infty}\frac{\mu_{Qp}}{\mu_p}>Q,$$
and $(\gamma_1)$ if
$$\hypertarget{gamma1}{(\gamma_1)}:\Leftrightarrow\sup_{p\in\NN_{>0}}\frac{\mu_p}{p}\sum_{k\ge p}\frac{1}{\mu_k}<\infty.$$
By \cite[Proposition 1.1]{petzsche} both conditions are equivalent for log-convex $M$ and for this proof condition \hyperlink{mnq}{$(\text{nq})$}, which is a general assumption in \cite{petzsche}, was not necessary, see also \cite[Theorem 3.11]{index}. In the literature \hyperlink{gamma1}{$(\gamma_1)$} is also called ''strong nonquasianalyticity condition.'' Moreover in \cite{Schmetsvaldivia00} the following generalization has been introduced (for $r\in\NN_{\ge 1}$):
$$\hypertarget{gammar}{(\gamma_r)}:\Leftrightarrow\sup_{p\in\NN_{>0}}\frac{(\mu_p)^{1/r}}{p}\sum_{k\ge p}\left(\frac{1}{\mu_k}\right)^{1/r}<\infty.$$
Of course, this condition makes sense for all $r>0$ and consequently $M$ has $(\gamma_r)$ if and only if $M^{1/r}$ has \hyperlink{gamma1}{$(\gamma_1)$}.

$(5)$ For two weight sequences $M=(M_p)_p$ and $N=(N_p)_p$ and $C>0$ we write $M\le CN$ if and only if $M_p\le CN_p$ holds for all $p\in\NN$. Moreover we define
$$M\hypertarget{mpreceq}{\preceq}N:\Leftrightarrow\;\exists\;C_1,C_2\ge 1\;\forall\;p\in\NN:\; M_p\le C_1 C_2^p N_p\Longleftrightarrow\sup_{p\in\NN_{>0}}\left(\frac{M_p}{N_p}\right)^{1/p}<\infty$$
and call them equivalent if
$$M\hypertarget{approx}{\approx}N:\Leftrightarrow\;M\hyperlink{mpreceq}{\preceq}N\;\text{and}\;N\hyperlink{mpreceq}{\preceq}M.$$
In the relations above one can replace $M$ and $N$ simultaneously by $m$ and $n$ because $M\hyperlink{mpreceq}{\preceq}N\Leftrightarrow m\hyperlink{mpreceq}{\preceq}n$.\vspace{6pt}

Let $d\in\NN_{>0}$ and $U\subseteq\RR^d$ be nonempty open, then for a weight sequence $M$ we introduce the (local) ultradifferentiable function class of Roumieu type by
\begin{equation*}\label{ultraRoumieu}
\mathcal{E}_{\{M\}}(U,\CC):=\{f\in\mathcal{E}(U,\CC):\;\forall\;K\subseteq U\;\text{compact}\;\exists\;h>0:\;\|f\|_{M,K,h}<\infty\},
\end{equation*}
and the Beurling type class by
\begin{equation*}\label{ultraBeurling}
\mathcal{E}_{(M)}(U,\CC):=\{f\in\mathcal{E}(U,\CC):\;\forall\;K\subseteq U\;\text{compact}\;\forall\;h>0:\;\|f\|_{M,K,h}<\infty\},
\end{equation*}
where we have put
\begin{equation*}\label{semi-norm-2}
\|f\|_{M,K,h}:=\sup_{k\in\NN^d,x\in K}\frac{|f^{(k)}(x)|}{h^{|k|} M_{|k|}}.
\end{equation*}
For compact sets $K$ with sufficiently regular boundary
$$\mathcal{E}_{M,h}(K,\CC):=\{f\in\mathcal{E}(K,\CC): \|f\|_{M,K,h}<\infty\}$$
is a Banach space where $\mathcal{E}(K,\CC)$ denotes the space of Whitney jets in $K$ which can be identified with the class of smooth functions on the interior $K^{\circ}$ with globally bounded derivatives. We have the topological vector space representations
\begin{equation*}\label{repr3}
\mathcal{E}_{\{M\}}(U,\CC)=\underset{K\subseteq U}{\varprojlim}\;\underset{h>0}{\varinjlim}\;\mathcal{E}_{M,h}(K,\CC)
\end{equation*}
and
\begin{equation*}\label{repr4}
\mathcal{E}_{(M)}(U,\CC)=\underset{K\subseteq U}{\varprojlim}\;\underset{h>0}{\varprojlim}\;\mathcal{E}_{M,h}(K,\CC).
\end{equation*}
From the definitions it is obvious that $M\hyperlink{mpreceq}{\preceq}N$ implies $\mathcal{E}_{[M]}(U,\CC)\subseteq\mathcal{E}_{[N]}(U,\CC)$.\vspace{6pt}

Let $M$ be a weight sequence and $U\subseteq\RR^d$ nonempty open
and $K$ a compact set with sufficiently regular boundary, then $\mathcal{E}_{[M]}(U,\CC)$ or $\mathcal{E}_{[M]}(K,\CC)$ are said to be {\itshape quasianalytic} if for all $x\in U$, or respectively for all $x\in K$, the map $\mathcal{B}_{x}$ is injective.\vspace{6pt}

Of course this definition can be extended to any subclass $\mathcal{S}\subseteq\mathcal{E}$. The following result characterizes nonquasianalyticity of ultradifferentiable classes in terms of the weight sequence $M$, e.g. see \cite[Proposition 4.3]{testfunctioncharacterization} where we have summarized the situation in a more general setting from \cite[Theorem 1.3.8]{hoermander}. For any weight sequence $M$, let $M^{\on{lc}}$ denote the log-convex minorant of $M$, e.g. see \cite[Def. 3.1, Prop. 3.2]{Komatsu73}.

\begin{proposition}\label{nonquasiremarks}
Let $M\in\RR_{>0}^{\NN}$ with $M_0=1$ and let $U\subseteq\RR^d$ be nonempty open. Then $\mathcal{E}_{[M]}(U,\CC)$ is nonquasianalytic if and only if $M^{\on{lc}}$ satisfies \hyperlink{mnq}{$(\on{nq})$} and in this case $\mathcal{C}^{\omega}(U,\CC)\subsetneq\mathcal{E}_{[M]}(U,\CC)=\mathcal{E}_{[M^{\on{lc}}]}(U,\CC)$ holds true.
\end{proposition}

\subsection{Relevant conditions for characterizing the surjectivity of the Borel map}\label{section32}
For $p,s\in\NN_{>0}$ given and two arbitrary sequences $M,N\in\RR_{>0}^{\NN}$ let
\begin{equation}\label{lambdaexpressionequ}
\lambda^{M,N}_{p,s}:=\sup_{0\le j<p}\left(\frac{M_p}{s^p N_j}\right)^{1/(p-j)}.
\end{equation}
We point out that the choice $j=0$ yields
\begin{equation}\label{lambdaexpressionequ0}
\lambda^{M,N}_{p,s}\ge\frac{1}{s(N_0)^{1/p}}(M_p)^{1/p},
\end{equation}
hence up to a constant a lower bound is always $(M_p)^{1/p}$. The next result proves a control from above (see \cite[$2(a)$]{surjectivity}).

\begin{lemma}\label{lambdaexpression}
Let $M$ and $N$ be log-convex weight sequences satisfying $M\le CN$ for some $C\ge 1$, then we get
\begin{equation}\label{lambdaexpressionequ1}
\forall\;p,s\in\NN_{>0}:\;\;\;\lambda^{M,N}_{p,s}\le C\min\{\mu_p,\nu_p\}.
\end{equation}
\end{lemma}

\demo{Proof}
Let $p,s\in\NN_{>0}$ and $0\le j\le p-1$. Then
\begin{align*}
\left(\frac{M_p}{s^p N_j}\right)^{1/(p-j)}&\le\left(\frac{M_p}{N_j}\right)^{1/(p-j)}\le\left(C\frac{M_p}{M_j}\right)^{1/(p-j)}=(C\mu_{j+1}\cdots\mu_p)^{1/(p-j)}\le C(\mu_p)^{(p-j)/(p-j)}=C\mu_p.
\end{align*}
Similarly we get
$$\left(\frac{M_p}{s^pN_j}\right)^{1/(p-j)}\le\left(C\frac{N_p}{N_j}\right)^{1/(p-j)}=(C\nu_{j+1}\cdots\nu_p)^{1/(p-j)}\le C(\nu_p)^{(p-j)/(p-j)}=C\nu_p.$$
\qed\enddemo

Let $M, N\in\hyperlink{LCset}{\mathcal{LC}}$ such that $M\le C N$ for some constant $C\ge 1$ and let $r\in\NN_{\ge 1}$. We consider now the following two conditions in the mixed weight sequence setting, for the definition even any $r>0$ makes sense:

\begin{equation*}\label{SVrami}
\hypertarget{SV}{(M,N)_{SV_r}}:\Longleftrightarrow\;\exists\;s\in\NN_{>0}:\;\;\sup_{p\in\NN_{>0}}\frac{(\lambda_{p,s}^{M,N})^{1/r}}{p}\sum_{k\ge p}\left(\frac{1}{\nu_k}\right)^{1/r}<\infty,
\end{equation*}
and
\begin{equation*}\label{gamma1rami}
\hypertarget{gammarmix}{(M,N)_{\gamma_r}}:\Longleftrightarrow\;\sup_{p\in\NN_{>0}}\frac{(\mu_p)^{1/r}}{p}\sum_{k\ge p}\left(\frac{1}{\nu_k}\right)^{1/r}<\infty.
\end{equation*}

So $(M,N)_{SV_r}$ is the $r$-ramification generalization of the characterizing condition $(\ast)$ in \cite[p. 385]{surjectivity} (taking $r=1$ yields $(\ast)$), whereas $(M,N)_{\gamma_r}$ is the generalization of condition \hyperlink{gammar}{$(\gamma_r)$} from \cite{Schmetsvaldivia00} to a mixed setting, since $(M,M)_{\gamma_r}$ is \hyperlink{gammar}{$(\gamma_r)$} for $M$.\vspace{6pt}

By Lemma \ref{lambdaexpression} we immediately get that $(M,N)_{\gamma_r}$ implies $(M,N)_{SV_r}$ for any $s\in\NN_{>0}$. If e.g. $M\in\hyperlink{LCset}{\mathcal{LC}}$ satisfies \hyperlink{mg}{$(\text{mg})$}, which is in this case equivalent to
\begin{equation}\label{mgequiv}
\exists\;A\ge 1\;\forall\;p\in\NN_{>0}:\;\;\;\mu_p\le A(M_p)^{1/p},
\end{equation}
e.g. see \cite[Lemma 2.2]{whitneyextensionweightmatrix} and the references therein, then by \eqref{lambdaexpressionequ0} and \eqref{lambdaexpressionequ1} conditions \hyperlink{SV}{$(M,N)_{SV_r}$} and \hyperlink{gammarmix}{$(M,N)_{\gamma_r}$} are equivalent.

\subsection{The r-interpolating sequence $P^{M,r}$}\label{rinterpolatingsequence}
Given a weight sequence $M\in\RR_{>0}^{\NN}$, in \cite[Lemma 2.3]{Schmetsvaldivia00} for any $r\in\NN_{\ge 1}$ the so-called $r${\itshape -interpolating sequence} $P^{M,r}$ was introduced as follows:
\begin{equation}\label{interpolatingsequ}
P^{M,r}_{rk+j}:=((M_k)^{r-j}(M_{k+1})^j)^{1/r},\hspace{20pt}\forall\;k\in\NN,\;\forall\;j\in\{0,\dots,r\}.
\end{equation}
We summarize some elementary facts: If $M_k\le C^kN_k$ for some $C\ge 1$ and all $k\in\NN$, then $P^{M,r}_{rk+j}\le(C^{1/r})^{rk+j}P^{N,r}_{rk+j}$. We have $P^{M,r}_{rn}=M_n$ for all $n\in\NN$ (i.e. for $r=1$ we get $P^{M,1}\equiv M$) and by denoting $\pi^{M,r}_k:=\frac{P^{M,r}_k}{P^{M,r}_{k-1}}$ (with $\pi^{M,r}_0:=1$) we see
\begin{equation}\label{interpolatingsequ0}
\forall\;k\in\NN\;\forall\;j\in\{1,\dots, r\}:\;\;\;\pi^{M,r}_{rk+j}=\frac{P^{M,r}_{rk+j}}{P^{M,r}_{rk+j-1}}=\left(\frac{M_k^{r-j}M_{k+1}^j}{M_k^{r-j+1}M_{k+1}^{j-1}}\right)^{1/r}=\left(\frac{M_{k+1}}{M_k}\right)^{1/r}=(\mu_{k+1})^{1/r}.
\end{equation}
Hence $M\in\hyperlink{LCset}{\mathcal{LC}}$ if and only if $P^{M,r}\in\hyperlink{LCset}{\mathcal{LC}}$ and moreover we can show:

\begin{lemma}\label{prmg}
Let $M\in\hyperlink{LCset}{\mathcal{LC}}$, $r\in\NN_{\ge 1}$ and $P^{M,r}$ be the $r$-interpolating sequence, then $M$ does have \hyperlink{mg}{$(\on{mg})$} if and only if $P^{M,r}$ does have \hyperlink{mg}{$(\on{mg})$}.
\end{lemma}

\demo{Proof}
It is well-known (e.g. see \cite[Lemma 2.2]{whitneyextensionweightmatrix}) that for $M\in\hyperlink{LCset}{\mathcal{LC}}$ condition \hyperlink{mg}{$(\on{mg})$} is equivalent to having $\sup_{p\ge 1}\frac{\mu_{2p}}{\mu_p}<\infty$.

On the one hand, for all $k\in\NN$ and $j\in\{1,\dots, r\}$ we have $2(rk+j)\le r(2k+1)+r$, hence by \eqref{interpolatingsequ0} we have $\frac{\pi^{M,r}_{2(rk+j)}}{\pi^{M,r}_{rk+j}}\le\frac{\pi^{M,r}_{r(2k+1)+r}}{\pi^{M,r}_{rk+j}}=\left(\frac{\mu_{2(k+1)}}{\mu_{k+1}}\right)^{1/r}$.

On the other hand, the choice $j=r$ in \eqref{interpolatingsequ0} yields $\frac{\mu_{2(k+1)}}{\mu_{k+1}}=\left(\frac{\pi^{M,r}_{2r(k+1)}}{\pi^{M,r}_{r(k+1)}}\right)^r$ for all $k\in\NN$.
\qed\enddemo

For arbitrary $r\in\NN_{\ge 1}$ by \eqref{interpolatingsequ0} we have that
\begin{equation*}\label{interpolatingsequnq1}
\sum_{l\ge 1}\frac{1}{\pi^{M,r}_l}=\sum_{k=0}^{\infty}\sum_{j=1}^r\frac{1}{\pi^{M,r}_{rk+j}}=\sum_{k=0}^{\infty}r\left(\frac{1}{\mu_{k+1}}\right)^{1/r},
\end{equation*}
hence \hyperlink{mnqr}{$(\text{nq}_r)$} holds for $M$ if and only if $P^{M,r}$ has \hyperlink{mnq}{($\text{nq})$}. So we immediately get by the classical Denjoy-Carleman theorem for ultradifferentiable functions (see Proposition \ref{nonquasiremarks}):

\begin{lemma}\label{prnonquasi}
Let $M\in\hyperlink{LCset}{\mathcal{LC}}$, $r\in\NN_{\ge 1}$ and $P^{M,r}$ be the $r$-interpolating sequence, then the following are equivalent:
\begin{itemize}
\item[$(i)$] $M$ satisfies \hyperlink{mnqr}{$(\on{nq}_r)$},

\item[$(ii)$] $P^{M,r}$ satisfies \hyperlink{mnq}{$(\on{nq})$},

\item[$(iii)$] $\mathcal{E}_{[P^{M,r}]}(U,\CC)$ is nonquasianalytic for every nonempty open set $U\subseteq\RR^d$.

\item[$(iv)$] $\mathcal{E}_{[P^{M,r}]}(K,\CC)$ is nonquasianalytic for every compact set $K\subseteq\RR^d$ with regular boundary.
\end{itemize}
\end{lemma}

Finally the next result generalizes \cite[Lemma 2.3 $(a)$]{Schmetsvaldivia00} to a mixed setting.

\begin{lemma}\label{prstrongnonquasi}
Let $M,N\in\hyperlink{LCset}{\mathcal{LC}}$, $r\in\NN_{\ge 1}$ and $P^{M,r}$, $P^{N,r}$ be the corresponding $r$-interpolating sequences. Then the following are equivalent:
\begin{itemize}
\item[$(i)$] \hyperlink{gammarmix}{$(M,N)_{\gamma_r}$},

\item[$(ii)$] \hyperlink{gammarmix}{$(P^{M,r},P^{N,r})_{\gamma_1}$}.
\end{itemize}
\end{lemma}

\demo{Proof}
We follow the proof given in \cite[Lemma 2.3 $(a)$]{Schmetsvaldivia00}.

$(i)\Rightarrow(ii)$ By using \eqref{interpolatingsequ0} we have
\begin{align*}
\sup_{l\in\NN_{\ge 1}}\frac{\pi^{M,r}_l}{l}\sum_{k\ge l}\frac{1}{\pi^{N,r}_k}&=\sup_{l\in\NN, j\in\{1,\dots,r\}}\frac{\pi^{M,r}_{rl+j}}{rl+j}\sum_{k\ge rl+j}\frac{1}{\pi^{N,r}_k}\le\sup_{l\in\NN}\frac{(\mu_{l+1})^{1/r}}{rl+1}\sum_{k\ge l}\sum_{j=1}^r\frac{1}{\pi^{N,r}_{rk+j}}
\\&
\le r\sup_{l\in\NN}\frac{(\mu_{l+1})^{1/r}}{rl+1}\sum_{k\ge l+1}\frac{1}{(\nu_k)^{1/r}}\le(1+r)\sup_{l\in\NN_{\ge 1}}\frac{(\mu_{l})^{1/r}}{l}\sum_{k\ge l}\frac{1}{(\nu_k)^{1/r}},
\end{align*}
where in the last estimate we have used that $\frac{r}{r(l-1)+1}\le\frac{1+r}{l}\Leftrightarrow r^2\le r^2l+1$ for all $l\in\NN_{\ge 1}$ and $r>0$.

$(ii)\Rightarrow(i)$ For all $n\in\NN_{\ge 1}$ we have
$$\frac{(\mu_{n})^{1/r}}{n}\sum_{k\ge n}\frac{1}{(\nu_k)^{1/r}}\le r\frac{\pi^{M,r}_{rn}}{rn}\sum_{k\ge n}\frac{1}{\pi^{N,r}_{rk}}\le r\frac{\pi^{M,r}_{rn}}{rn}\sum_{k\ge rn}\frac{1}{\pi^{N,r}_{k}},$$
taking into account again \eqref{interpolatingsequ0}.
\qed\enddemo

\subsection{Associated sequence and test function spaces}\label{section33}
Let $M\in\RR_{>0}^{\NN}$ be a weight sequence and for a given sequence $\mathbf{a}:=(a_p)_p\in\CC^{\NN}$ we put
$$|\mathbf{a}|_{M,h}:=\sup_{p\in\NN}\frac{|a_p|}{h^p M_p},$$
and introduce $\Lambda_{M,h}:=\{(a_p)_p\in\CC: |\mathbf{a}|_{M,h}<\infty\}$. Furthermore we set
$$\Lambda_{(M)}:=\{(a_p)_p\in\CC^{\NN}: \forall\;h>0: |\mathbf{a}|_{M,h}<\infty\},$$
and
$$\Lambda_{\{M\}}:=\{(a_p)_p\in\CC: \exists\;h>0: |\mathbf{a}|_{M,h}<\infty\}.$$
We endow $\Lambda_{(M)}$ resp. $\Lambda_{\{\mathcal{M}\}}$ with a natural projective, respectively inductive, topology via $$\Lambda_{(M)}=\underset{h>0}{\varprojlim}\;\Lambda_{M,h},\hspace{30pt}\Lambda_{\{\mathcal{M}\}}=\underset{h>0}{\varinjlim}\;\Lambda_{M,h}.$$
If $M\in\hyperlink{LCset}{\mathcal{LC}}$, then both spaces are rings with respect to convolution $(\mathbf{a}\star \mathbf{b})_n:=\sum_{k=0}^n a_k b_{n-k}$:
\begin{align*}
|(\mathbf{a}\star \mathbf{b})|_{M,2h}&=\sup_n\frac{|(\mathbf{a}\star\mathbf{b})_n|}{(2h)^n M_n}\le\sup_n\frac{\sum_{k=0}^n|a_k||b_{n-k}|}{(2h)^n M_n}\le\sup_n\frac{1}{2^n}\sum_{k=0}^n\biggl(\underbrace{\frac{|a_k|}{h^k M_k}}_{\le|\mathbf{a}|_{M,h}}\underbrace{\frac{|b_{n-k}|}{h^{n-k} M_{n-k}}}_{\le|\mathbf{b}|_{M,h}}\biggr)
\\&
\le\sup_n\frac{n+1}{2^n}|\mathbf{a}|_{M,h}|\mathbf{b}|_{M,h}\le|\mathbf{a}|_{M,h}|\mathbf{b}|_{M,h}<\infty.
\end{align*}
In this estimate we have used that $M$ is log-convex and normalized, so $M_k M_{n-k}\le M_n$ for all $n,k\in\NN$ with $k\le n$. From the definitions it is obvious that $M\hyperlink{mpreceq}{\preceq}N$ implies $\Lambda_{[M]}\subseteq\Lambda_{[N]}$.\vspace{6pt}

Let $M$ be a weight sequence and $r\in\NN_{\ge 1}$. Then for each $h>0$ and $a>0$ we define the Banach space
\begin{align*}
&\mathcal{D}_{r,M,h}([-a,a]):=
\\&
\{f\in\mathcal{E}(\RR,\CC):\;\supp(f)\subseteq[-a,a],\;f^{(rn+j)}(0)=0,\;\forall\;n\in\NN\;\forall\;j\in\{1,\dots,r-1\}, \sup_{n\in\NN, x\in\RR}\frac{|f^{(rn)}(x)|}{h^nM_n}<\infty\},
\end{align*}
and the Roumieu type class
\begin{equation}\label{Roumieu-LB-space}
\mathcal{D}_{r,\{M\}}([-a,a]):=\underset{h>0}{\varinjlim}\;\mathcal{D}_{r,M,h}([-a,a]),
\end{equation}
which is a countable $(LB)$-space, respectively
\begin{equation}\label{Beurling-Frechet-space}
\mathcal{D}_{r,(M)}([-a,a]):=\underset{h>0}{\varprojlim}\;\mathcal{D}_{r,M,h}([-a,a]),
\end{equation}
which is a Frechét space (see \cite[Section 3]{Schmetsvaldivia00}). For convenience we put $\|f\|_{r,M,h}:=\sup_{n\in\NN, x\in\RR}\frac{|f^{(rn)}(x)|}{h^nM_n}$. If $r=1$, then we precisely obtain the spaces considered in \cite{surjectivity} and $M\hyperlink{mpreceq}{\preceq}N$ implies $\mathcal{D}_{r,[M]}([-a,a])\subseteq\mathcal{D}_{r,[N]}([-a,a])$.

\begin{remark}\label{remarktestfunctionnontrivial}
We finish this section with the following comments:
\begin{itemize}
\item[$(i)$] In the main results Theorem \ref{Roumieu-Surjectivitytheorem} and Theorem \ref{Beurling-Surjectivitytheorem} equivalently we could replace $\mathcal{D}_{r,[M]}([-1,1])$ by $\mathcal{D}_{r,[M]}([-a,a])$ for arbitrary $a>0$: The spaces obviously do coincide as sets by composing the functions with a dilation, see also Lemma \ref{Petzsche2remark} below.

\item[$(ii)$] It is not automatically clear that the classes introduced in \eqref{Roumieu-LB-space} and \eqref{Beurling-Frechet-space} are nontrivial, i.e. $\neq\{0\}$. In Section \ref{nonquasi} we will characterize the nontriviality in terms of $M$ as in the classical Denjoy-Carleman theorem and we will see that this question is characterized by the nonquasianalyticity of $P^{M,r}$, see Lemma \ref{prnonquasi} and also Lemma \ref{Petzsche2remark} below.
\end{itemize}
\end{remark}

\section{The image of the Borel mapping in the Roumieu case}\label{section4}
Let, from now on in this section, $M,N\in\RR_{>0}^{\NN}$ and $r\in\NN_{\ge 1}$ be such that

\begin{itemize}
\item[$(I)$] $M,N\in\hyperlink{LCset}{\mathcal{LC}}$,


\item[$(II)_{R,r}$] $\liminf_{p\rightarrow\infty}\left(\frac{(M_p)^{1/r}}{p!}\right)^{1/p}>0$ (the letter $R$ in the notation stands for {\itshape Roumieu}),

\item[$(III)$] $M\hyperlink{mpreceq}{\preceq}N$, i.e. $M_p\le C^p N_p$ for a constant $C\ge 1$ and all $p\in\NN$.
\end{itemize}

\vspace{6pt}
Concerning these conditions we give the following comments.\vspace{6pt}

\begin{remark}\label{roumieucondtionremark}
\begin{itemize}
\item[$(i)$] Replacing $N=(N_p)_p$ by $L=(L_p)_p$ with $L_p:=C^pN_p$, $C$ denoting the constant from $(III)$, we have $M_p\le L_p$ for all $p\in\NN$. Of course $L\hyperlink{approx}{\approx}N$ and so we can assume from now on without loss of generality that even $M\le N$ holds true (which will simplify the notation).

\item[$(ii)$] For all $M\in\hyperlink{LCset}{\mathcal{LC}}$ and $r>1$ condition $(II)_{R,r}$ does imply $(II)_{R,1}$, indeed $\lim_{p\rightarrow\infty}(m_p)^{1/p}=\infty$ since $(M_p)^{1/p}\rightarrow\infty$ as $p\rightarrow\infty$ and $m_p=\frac{(M_p)^{1/r}}{p!}(M_p)^{1-\frac{1}{r}}$.

\item[$(iii)$] $(II)_{R,r}$ for $M$ does imply $(II)_{R,1}$ for $P^{M,r}$:

If $r=1$, then $M\equiv P^{M,1}$ and nothing is to prove. If $r\in\NN_{\ge 2}$, then we have $M_p\ge C^{rp}p!^r$ for some $0<C\le 1$ and all $p\in\NN$, hence by \eqref{interpolatingsequ}
\begin{align*}
P^{M,r}_{rk+j}\ge(C^{(kr)(r-j)}k!^{r(r-j)}C^{rj(k+1)}(k+1)!^{rj})^{1/r}=C^{kr+j}(k!)^{r-j}(k+1)!^j
\end{align*}
for all $k\in\NN$ and $1\le j\le r$. Then
$$(rk+j)!\le 2^{rk+j}j!(rk)!\le 2^{rk+j}r!(rk)!\le 2^{rk+j}r!C_1^{rk}k!^r\le r!(2C_1)^{rk+j}k!^{r-j}(k+1)!^j,$$
where $C_1$ is a suitable positive constant depending only on $r$; consequently
$$P^{M,r}_{rk+j}\ge r!^{-1}(C/(2C_1))^{rk+j}(rk+j)!.$$
Thus $\liminf_{k\rightarrow\infty}(P^{M,r}_k/k!)^{1/k}>0$ is shown.
\end{itemize}
\end{remark}

The goal is to prove the following characterization which is generalizing \cite[Theorem 1.1]{surjectivity}.

\begin{theorem}\label{Roumieu-Surjectivitytheorem}
Let $M$ and $N$ be as assumed above and $r\in\NN_{\ge 1}$. Then
\begin{equation}\label{Roumieu-Surjectivitytheoreminclusion}
\Lambda_{\{M\}}\subseteq\{(f^{(rj)}(0))_{j\in\NN}: f\in\mathcal{D}_{r,\{N\}}([-1,1])\}
\end{equation}
if and only if \hyperlink{SV}{$(M,N)_{SV_r}$}.
\end{theorem}

\begin{remark}\label{Roumieu-Surjectivitytheoremremark}
This theorem extends \cite[Theorem 1.1]{surjectivity} also to general $r$-interpolating spaces because there only the case $r=1$ was considered. But even in this case our approach is slightly stronger than the result from \cite{surjectivity} since our assumptions on $M$ and $N$ are more general. More precisely:

\begin{itemize}
\item[$(i)$] We only need $\liminf_{k\rightarrow\infty}(m_k)^{1/k}>0$ for $M$ instead of $\lim_{k\rightarrow\infty}(m_k)^{1/k}=\infty$, which is the general assumption $(d)$ on \cite[p. 385]{surjectivity},

\item[$(ii)$] We only require $M\hyperlink{mpreceq}{\preceq}N$ instead of the stronger assumption $\mu_p\le\nu_p$ for all $p\in\NN$ (which is assumption $(c)$ on \cite[p. 385]{surjectivity}).

\item[$(iii)$] Most importantly, it is not required to have condition \hyperlink{mnq}{$(\on{nq})$} for $N$. We give here the argument for arbitrary $r\in\NN_{\ge 1}$: If $M$ and $N$ are related by \hyperlink{SV}{$(M,N)_{SV_r}$}, then by choosing $p=1$ there (which yields $(\lambda^{M,N}_{1,s})^{1/r}=\left(\frac{M_1}{s N_0}\right)^{1/r}$) we immediately get that $N$ has to satisfy \hyperlink{mnqr}{$(\on{nq}_r)$}. Conversely, if \eqref{Roumieu-Surjectivitytheoreminclusion} holds true, then it is obvious that the class $\mathcal{D}_{r,\{N\}}([-1,1])$ has to be nontrivial because $\Lambda_{\{M\}}$ does contain (at least) all sequences with finitely many entries $\neq 0$. Consequently, by Theorem \ref{nonquasitheorem} below, this nontriviality of $\mathcal{D}_{r,\{N\}}([-1,1])$ does imply condition \hyperlink{mnqr}{$(\on{nq}_r)$} for $N$.
\end{itemize}
\end{remark}

First we recall \cite[Theorem 2.2]{petzsche}, which follows from \cite[Theorem 1.3.5]{hoermander}, for a proof see also \cite[Lemma 5.1.6]{diploma}.

\begin{lemma}\label{Petzsche2}
Let $M\in\hyperlink{LCset}{\mathcal{LC}}$ and assume that $a:=\sum_{j=1}^{\infty}\frac{1}{\mu_j}<\infty$. Then there exists a smooth function $\varphi$ whose support is contained in $[-a,a]$, such that $0\le\varphi(x)\le 1$ for all $x\in[-a,a]$, and $\varphi^{(j)}(0)=\delta_{j,0}$ (Kronecker delta). Furthermore we have $\big\|\varphi^{(j)}\big\|_{\infty}\le 2^jM_j$ for all $j\in\NN$.\vspace{12pt}

In particular one can say: $\varphi$ is a nontrivial function ($\varphi(0)=1$) with compact support and $\varphi\in\mathcal{E}_{\{M\}}(\RR,\CC)$ (take $h=2$). Thus the ultradifferentiable class $\mathcal{E}_{\{M\}}$ is nonquasianalytic.
\end{lemma}

In the next statement we will use the previous result to justify that the class $\mathcal{D}_{r,\{M\}}([-a,a])$ defined in \eqref{Roumieu-LB-space} is nontrivial. As mentioned before, for a complete characterization (also for the Beurling case) we refer to Section \ref{nonquasi} below and we will have to make use of the following construction (for the Roumieu case) in the first main result Theorem \ref{Roumieu-Surjectivitytheorem1}.

\begin{lemma}\label{Petzsche2remark}
Let $M\in\hyperlink{LCset}{\mathcal{LC}}$ and $r\in\NN_{\ge 1}$ be given and assume that \hyperlink{mnqr}{$(\on{nq}_r)$} holds true. Then we get for all $b>0$ that $\mathcal{D}_{r,[M]}([-b,b])$ is nontrivial.
\end{lemma}

\demo{Proof}
We set $a:=\sum_{k=0}^{\infty}\left(\frac{1}{\mu_{k+1}}\right)^{1/r}<\infty$ and apply Lemma \ref{Petzsche2} to the $r$-interpolating sequence $P^{M,r}$ from \eqref{interpolatingsequ} (see Lemma \ref{prnonquasi}). Thus $\varphi$ belongs to the class $\mathcal{D}_{\{P^{M,r}\}}([-a,a])$ and even to $\mathcal{D}_{r,\{M\}}([-a,a])$ (since $\varphi^{(j)}(0)=\delta_{j,0}$ and $P^{M,r}_{rj}=M_j$ for all $j\in\NN$) and has compact support in $[-a,a]$. Finally we get $\big\|\varphi^{(rj)}\big\|_{\infty}\le 2^{rj}P^{M,r}_{rj}=2^{rj}M_j$ for all $j\in\NN$.\vspace{6pt}

By making a rescaling/dilation we can achieve that $\varphi$ has support contained in $[-b,b]$ for any $b>0$.


For the Beurling type we have to recall that in the proof of \cite[Theorem 2.1 $(a)(i)$]{petzsche} even a sequence $(\chi_p)_p$ of functions with compact support in the ultradifferentiable class $\mathcal{E}_{(P^{M,r})}(\RR,\CC)$ has been constructed which satisfies $\chi^{(j)}_p(0)=\delta_{j,p}$. So the corresponding result holds true for the Beurling type classes $\mathcal{D}_{r,(M)}([-b,b])$ as well (by choosing $\chi_0$ or more general some $\chi_{rp}$) and making again a rescaling.
\qed\enddemo

Using this preparation we are able to generalize \cite[Theorem 3.2]{surjectivity}. We are going to follow the original proof and make adjustments where necessary.

\begin{theorem}\label{Roumieu-Surjectivitytheorem1}
Let $M$ and $N$ be as assumed above and $r\in\NN_{\ge 1}$ be given. If \hyperlink{SV}{$(M,N)_{SV_r}$} holds true, then there exists $d>0$ such that for all $c\in\NN_{>0}$ there exists a continuous linear extension map $T_c:\Lambda_{M,c}\rightarrow\mathcal{D}_{r,N,cd}([-1,1])$, which implies
$$\Lambda_{\{M\}}\subseteq\{(f^{(rj)}(0))_{j\in\NN}: f\in\mathcal{D}_{r,\{N\}}([-1,1])\}.$$
\end{theorem}

\demo{Proof}
For convenience we write $\lambda_{p,s}$ instead of $\lambda^{M,N}_{p,s}$. Let $h\ge 1$ (large) be arbitrary but fixed and $s\in\NN_{>0}$ be coming from \hyperlink{SV}{$(M,N)_{SV_r}$}. For $p\in\NN$ we consider the sequence $\tau^p:=(\tau^p_j)_{j\ge 0}$ defined by
$$\tau^p:=\big(1,\underbrace{(h\lambda_{p,s})^{1/r},\dots,(h\lambda_{p,s})^{1/r}}_{2pr-\text{times}},\underbrace{(h\nu_{2p+1})^{1/r},\dots,(h\nu_{2p+1})^{1/r}}_{r-\text{times}},\underbrace{(h\nu_{2p+2})^{1/r},\dots,(h\nu_{2p+2})^{1/r}}_{r-\text{times}},\dots\big).$$
By \eqref{lambdaexpressionequ1} (with $C=1$) we see that each $\tau^p$ is increasing since $\lambda_{p,s}\le\nu_p\le\nu_{2p+1}$ and so $T^p=(T^p_j)_j$ defined by $T^p_j:=\prod_{i=0}^j\tau^p_i$ is log-convex (and its quotients are tending to infinity).\vspace{6pt}

{\itshape The case $p=0$.} Since $B:=h^{1/r}\sum_{j=1}^{\infty}\left(\frac{1}{h\nu_j}\right)^{1/r}=\sum_{j=1}^{\infty}\left(\frac{1}{\nu_j}\right)^{1/r}<\infty$ we can apply Lemma \ref{Petzsche2remark} for $a=B$, $b=B/h^{1/r}$ and so $c=h^{1/r}$. So there exists $\varrho_{h,0}\in\mathcal{D}_{r,\{N\}}([-B/(h^{1/r}),B/(h^{1/r})])$ with $0\le\varrho_{h,0}\le 1$, $\varrho_{h,0}^{(j)}(0)=\delta_{j,0}$ for all $j\in\NN$ and finally $|\varrho_{h,0}^{(rj)}(t)|\le 2^{rj}h^jN_j$ for all $j\in\NN$ and $t\in\RR$.\vspace{6pt}

{\itshape The case $p\in\NN_{>0}$.} Property \hyperlink{SV}{$(M,N)_{SV_r}$} means that there exists some $A>1$ such that for all $p\in\NN_{>0}$ we can estimate as follows
\begin{align*}
&\sum_{k=1}^{\infty}\frac{1}{\tau^p_k}=\frac{2pr}{(h\lambda_{p,s})^{1/r}}+r\sum_{k=2p+1}^{\infty}\left(\frac{1}{h\nu_k}\right)^{1/r}=\frac{2pr}{(h\lambda_{p,s})^{1/r}}\left(1+\frac{(\lambda_{p,s})^{1/r}}{2p}\sum_{k=2p+1}^{\infty}\left(\frac{1}{\nu_k}\right)^{1/r}\right)
\\&
\le\frac{2pr}{(h\lambda_{p,s})^{1/r}}\left(1+\frac{(\lambda_{p,s})^{1/r}}{2p}\sum_{k=p}^{\infty}\left(\frac{1}{\nu_k}\right)^{1/r}\right)\le\frac{2rp}{(h\lambda_{p,s})^{1/r}}A.
\end{align*}
Applying Lemma \ref{Petzsche2} for each sequence $T^p$, $p\in\NN_{>0}$, we get that there exists

$\varrho_{h,p}\in\mathcal{D}([-2Arp/(h\lambda_{p,s})^{1/r},2Arp/(h\lambda_{p,s})^{1/r}])$ with $0\le\varrho_{h,p}\le 1$, $\varrho_{h,p}^{(j)}(0)=\delta_{j,0}$ for all $j\in\NN$ and satisfying for all $t\in\RR$:
\begin{equation}\label{varrhoinequ}
|\varrho_{h,p}^{(j)}(t)|\le
\begin{cases}
2^j(h\lambda_{p,s})^{j/r} &\text{if}\;1\le j\le 2rp,
\\
2^jh^{j/r}(\lambda_{p,s})^{2p}\prod_{k=1}^{k_j}\nu_{2p+k}(\nu_{2p+k_j+1})^{(j-r(2p+k_j))/r} &\text{if}\;j>2rp,
\end{cases}
\end{equation}
where for the second inequality we have put $k_j\in\NN$ satisfying $(2p+k_j)r<j\le(2p+k_j+1)r$ (if $k_j=0$, then the product is understood to be equal $1$, so $|\varrho_{h,p}^{(j)}(t)|\le 2^jh^{j/r}(\lambda_{p,s})^{2p}(\nu_{2p+1})^{(j-2pr)/r}$). We introduce the smooth function $\chi:\RR\rightarrow\RR$ defined by
\begin{equation}\label{functionchi}
\chi_{h,rp}(t):=\varrho_{h,p}(t)\frac{t^{rp}}{(rp)!},
\end{equation}
and in the next step we have to estimate all $(rj)$-derivatives of $\chi_{h,pr}$ for each $p\in\NN$.

{\itshape The case $p=0$.} We have $|\chi_{h,0}^{(rj)}(t)|=|\varrho_{h,0}^{(rj)}(t)|\le 2^{rj}h^jN_j=2^{rj}h^j\frac{N_j}{M_0}$ for all $t\in\RR$ and $j\in\NN$.

{\itshape The case $p\in\NN_{>0}$:} We are going to prove
\begin{equation}\label{equ1}
\forall\;t\in\RR\;\forall\;j\in\NN:\;\;|\chi_{h,rp}^{(rj)}(t)|\le \frac{N_j}{M_p}\left(\frac{2A\exp(1)s}{h^{1/r}}\right)^{rp}h^j\left(2+\frac{1}{2A}\right)^{rj}.
\end{equation}
First the Leibniz-formula gives (since $(t^p)^{(j)}=0$ for any $j>p$) that:
\begin{equation}\label{leibniz}
|\chi_{h,rp}^{(rj)}(t)|\le\left(\frac{2A\exp(1)}{h^{1/r}}\right)^{rp}\sum_{l=\max\{0,r(j-p)\}}^{rj}\binom{rj}{l}|\varrho_{h,p}^{(l)}(t)|\left(\frac{2A}{h^{1/r}}\right)^{l-rj}\left(\frac{1}{(\lambda_{p,s})^{1/r}}\right)^{rp+l-rj}.
\end{equation}
We point out that $t\in[-2Arp/(h\lambda_{p,s})^{1/r},2Arp/(h\lambda_{p,s})^{1/r}]$ and moreover $rj-rp\le l$ and $l\le rj$ imply $0\le rp+l-rj\le rp$. So $\frac{(rp)^{rp+l-rj}}{(rp+l-rj)!}\le\frac{(rp)^{rp}}{(rp)!}\le\exp(rp)$ because $\frac{p^m}{m!}\le\frac{p^{m+1}}{(m+1)!}\Leftrightarrow m+1\le p$ and $\frac{p^p}{p!}\le\exp(p)=\sum_{m=0}^{\infty}\frac{p^m}{m!}$ for all $p\in\NN$.\vspace{6pt}

Now we have to distinguish between two cases:\vspace{6pt}

{\itshape Case 1, $p\in\NN_{>0}$ and $0\le j\le 2p\Leftrightarrow 0\le rj\le 2rp$.} By \eqref{varrhoinequ} we have $|\varrho_{h,p}^{(l)}(t)|\le 2^l(h\lambda_{p,s})^{l/r}$ for all $l\in\NN$ with $\max\{0,r(j-p)\}\le l\le jr$ (if $l=0$, then this holds true by having $\|\varrho_{h,p}\|_{\infty}\le 1$) and so
\begin{align*}
|\chi_{h,rp}^{(rj)}(t)|&\le\left(\frac{2A\exp(1)}{h^{1/r}}\right)^{rp}h^j(\lambda_{p,s})^{j-p}\sum_{l=\max\{0,r(j-p)\}}^{rj}\binom{rj}{l}2^l(2A)^{l-rj}
\\&
\le\left(\frac{2A\exp(1)}{h^{1/r}}\right)^{rp}h^j(\lambda_{p,s})^{j-p}\left(2+\frac{1}{2A}\right)^{rj}.
\end{align*}
{\itshape Case 1, subcase $(a)$, $0\le j<p$.} By definition \eqref{lambdaexpressionequ} we have $(\lambda_{p,s})^{p-j}\ge\frac{M_p}{s^pN_j}$ and so the previous estimate shows (by having $s\in\NN_{\ge 1}$)
\begin{equation*}\label{equ2}
|\chi_{h,rp}^{(rj)}(t)|\le\frac{N_j}{M_p}\left(\frac{2A\exp(1)s}{h^{1/r}}\right)^{rp}h^j\left(2+\frac{1}{2A}\right)^{rj}\hspace{20pt}\text{for}\;0\le j<p.
\end{equation*}
{\itshape Case 1, subcase $(b)$, $j=p$.} In this case $(\lambda_{p,s})^{p-j}=1$ and since $\frac{N_j}{M_p}=\frac{N_p}{M_p}\ge 1$ we have
\begin{equation*}\label{equ3}
|\chi_{h,rp}^{(rj)}(t)|\le \frac{N_j}{M_p}\left(\frac{2A\exp(1)}{h^{1/r}}\right)^{rp}h^j\left(2+\frac{1}{2A}\right)^{rj}\hspace{20pt}\text{for}\;j=p.
\end{equation*}
{\itshape Case 1, subcase $(c)$, $p<j\le 2p$.} First, by \eqref{lambdaexpressionequ1} we get $\lambda_{p,s}\le\nu_p$ for all $p,s\in\NN_{>0}$. Since $j-p>0$ we have
\begin{equation}\label{subcase1equ}
\forall\;j,p\in\NN_{>0},\;j>p:\;\;\;(\lambda_{p,s})^{j-p}\le(\nu_p)^{j-p}\le\nu_{p+1}\cdots\nu_j=\frac{N_j}{N_p}\le\frac{N_j}{M_p}.
\end{equation}
This proves
\begin{equation*}\label{equ4}
|\chi_{h,rp}^{(rj)}(t)|\le\frac{N_j}{M_p}\left(\frac{2A\exp(1)}{h^{1/r}}\right)^{rp}h^j\left(2+\frac{1}{2A}\right)^{rj}\hspace{20pt}\text{for}\;p<j\le 2p.
\end{equation*}

{\itshape Case 2, $p\in\NN_{>0}$ and $2p<j\Leftrightarrow 2rp<rj$.} By \eqref{varrhoinequ} we have

$|\varrho_{h,p}^{(l)}(t)|\le 2^lh^{l/r}(\lambda_{p,s})^{2p}\prod_{k=1}^{k_l}\nu_{2p+k}(\nu_{2p+k_l+1})^{(l-r(2p+k_l))/r}$ for all $l>2rp$ and $k_l\in\NN$ satisfying $(2p+k_l)r<l\le(2p+k_l+1)r$. In this case we are interested in such values satisfying $0\le k_l\le j-2p-1$. In the estimate we decompose the sum $\sum_{l=r(j-p)}^{rj}$ in \eqref{leibniz} into $\sum_{l=r(j-p)}^{2rp}+\sum_{l=2rp+1}^{rj}$. Hence
\begin{align*}
&|\chi_{h,rp}^{(rj)}(t)|\le\left(\frac{2A\exp(1)}{h^{1/r}}\right)^{rp}\sum_{l=r(j-p)}^{2rp}\binom{rj}{l}\underbrace{|\varrho_{h,p}^{(l)}(t)|}_{\le 2^l(h\lambda_{p,s})^{l/r}}\left(\frac{2A}{h^{1/r}}\right)^{l-rj}\left(\frac{1}{(\lambda_{p,s})^{1/r}}\right)^{rp+l-rj}
\\&
+\left(\frac{2A\exp(1)}{h^{1/r}}\right)^{rp}\sum_{l=2rp+1}^{rj}\binom{rj}{l}\underbrace{|\varrho_{h,p}^{(l)}(t)|}_{\le 2^lh^{l/r}(\lambda_{p,s})^{2p}\prod_{k=1}^{k_l}\nu_{2p+k}(\nu_{2p+k_l+1})^{(l-r(2p+k_l))/r}}\left(\frac{2A}{h^{1/r}}\right)^{l-rj}\left(\frac{1}{(\lambda_{p,s})^{1/r}}\right)^{rp+l-rj}
\\&
\le\left(\frac{2A\exp(1)}{h^{1/r}}\right)^{rp}h^j\underbrace{(\lambda_{p,s})^{j-p}}_{\le N_j/M_p}\sum_{l=r(j-p)}^{2rp}\binom{rj}{l}2^l(2A)^{l-rj}
\\&
+\left(\frac{2A\exp(1)}{h^{1/r}}\right)^{rp}h^j\sum_{l=2rp+1}^{rj}\binom{rj}{l}2^l(2A)^{l-rj}\underbrace{(\lambda_{p,s})^{p+j-l/r}\prod_{k=1}^{k_l}\nu_{2p+k}(\nu_{2p+k_l+1})^{(l-r(2p+k_l))/r}}_{\le N_j/M_p (\star)}
\\&
\le\frac{N_j}{M_p}\left(\frac{2A\exp(1)}{h^{1/r}}\right)^{rp}h^j\sum_{l=r(j-p)}^{rj}\binom{rj}{l}2^l(2A)^{l-rj}\le \frac{N_j}{M_p}\left(\frac{2A\exp(1)}{h^{1/r}}\right)^{rp}h^j\left(2+\frac{1}{2A}\right)^{rj}.
\end{align*}
In the first sum we have used again \eqref{subcase1equ} (note that $j>p$).

To prove $(\star)$ we have to estimate as follows: First, since $\nu_p\ge 1$ for all $p\in\NN$ and $p+j-l/r\ge 0\Leftrightarrow r(p+j)\ge l$ (since $l\le rj$ in the second sum), we get $(\lambda_{p,s})^{p+j-l/r}\le(\nu_p)^{p+j-l/r}$ (by \eqref{lambdaexpressionequ1}) and so $(\lambda_{p,s})^{p+j-l/r}\le\nu_{p+1}\cdots\nu_{2p}(\nu_p)^{j-l/r}$. Second we have by log-convexity and since $j-l/r\ge 0\Leftrightarrow rj\ge l$:
\begin{align*}
(\nu_{2p+k_l+1})^{(l-r(2p+k_l))/r}(\nu_p)^{j-l/r}&\le(\nu_{2p+k_l+1})^{(l-r(2p+k_l))/r}(\nu_{2p+k_l+1})^{j-l/r}=(\nu_{2p+k_l+1})^{j-2p-k_l}
\\&
\le\nu_{2p+k_l+1}\cdots\nu_j.
\end{align*}
Thus by combining both estimations we get as desired
\begin{align*}
(\lambda_{p,s})^{p+j-l/r}\prod_{k=1}^{k_l}\nu_{2p+k}(\nu_{2p+k_l+1})^{(l-r(2p+k_l))/r}\le\nu_{p+1}\cdots\nu_j=\frac{N_j}{N_p}\le\frac{N_j}{M_p}.
\end{align*}

This finishes the proof of \eqref{equ1}.\vspace{12pt}

The case $j=0$ and $N_0=1$ in \eqref{lambdaexpressionequ} yield $(M_p)^{1/(rp)}\le(s\lambda_{p,s})^{1/r}$ for all $p\in\NN_{>0}$. By assumption $(II)_{R,r}$ we have $\liminf_{p\rightarrow\infty}\left(\frac{(M_p)^{1/r}}{p!}\right)^{1/p}>0$ and by applying Stirling's formula $\liminf_{p\rightarrow\infty}\left(\frac{(M_p)^{1/r}}{p^p}\right)^{1/p}>0$, too. Thus $\frac{(M_p)^{1/(rp)}}{p}\le\frac{(s\lambda_{p,s})^{1/r}}{p}\Leftrightarrow\frac{p}{(\lambda_{p,s})^{1/r}}\le s^{1/r}\frac{p}{(M_p)^{1/(rp)}}$ and we get $\sup_{p\in\NN_{>0}}\frac{p}{(\lambda_{p,s})^{1/r}}<\infty$.

We can choose now some number $l\in\NN_{>0}$ large enough, depending on given $s$ and $r$ and satisfying
\begin{itemize}
\item[$(a)$] $\frac{B}{l}<\frac{1}{2}$,

\item[$(b)$] $\frac{2Arp}{l(\lambda_{p,s})^{1/r}}<1$ for all $p\in\NN_{>0}$ and

\item[$(c)$] $\frac{2A\exp(1)s}{l^{1/r}}<\frac{1}{2}$.
\end{itemize}
We may suppose from now on that $h\in\NN_{>0}$ and $h\ge l^r$. For this particular $h$ we summarize:
\begin{itemize}
\item[$(i)$] By $(a)$ and $(b)$ the support of $\varrho_{h,p}$ is contained in $[-1,1]$ for all $p\in\NN$ since $\frac{B}{h^{1/r}}\le\frac{B}{l}$ and $\frac{2Arp}{(h\lambda_{p,s})^{1/r}}\le\frac{2Arp}{l(\lambda_{p,s})^{1/r}}$.
\item[$(ii)$] As shown in \eqref{equ1} we have
$$\forall\;t\in\RR\;\forall\;j,p\in\NN:\;\;\;|\chi_{h,rp}^{(rj)}(t)|\le \frac{N_j}{M_p}\left(\frac{2A\exp(1)s}{h^{1/r}}\right)^{rp}h^j\left(2+\frac{1}{2A}\right)^{rj}.$$
\item[$(iii)$] $\varrho_{h,p}^{(l)}(0)=\delta_{l,0}$ for all $p\in\NN$ (and arbitrary $h$) and definition \eqref{functionchi} imply that $$\forall\;j,p\in\NN:\;\chi^{(j)}_{h,rp}(0)=\delta_{rp,j}.$$
\end{itemize}
$(i)-(iii)$ prove that $\chi_{h,rp}\in\mathcal{D}_{r,\{N\}}([-1,1])$ holds true (for all $p\in\NN$).

We put $d:=l^{1/r}(2+1/(2A))$ (also depending on $s$ and $r$ via chosen $l$) and let $c\in\NN_{>0}$ be arbitrary but fixed. Consider a sequence $\mathbf{a}=(a_p)_p\in\Lambda_{M,c}$, then
\begin{align*}
\forall\;p,j\in\NN\;\forall\;t\in\RR:\;\;|a_p\chi^{(rj)}_{cl,rp}(t)|&\le|\mathbf{a}|_{M,c}c^p M_p\frac{N_j}{M_p}(cl)^j\left(\frac{2A\exp(1)s}{(cl)^{1/r}}\right)^{rp}\left(2+\frac{1}{2A}\right)^{rj}
\\&
\underbrace{\le}_{(c)}|\mathbf{a}|_{M,c}c^pN_j(cl)^j\frac{1}{2^{rp}c^p}\left(2+\frac{1}{2A}\right)^{rj}\le\frac{1}{2^{rp}}|\mathbf{a}|_{M,c}c^jd^{rj}N_j,
\end{align*}
where the last inequality holds by the choice of $d$. Note that by the choice of $l$ in $(c)$ we have $\frac{2A\exp(1)s}{l^{1/r}}<\frac{1}{2}$ which is equivalent to $\left(\frac{2A\exp(1)s}{(cl)^{1/r}}\right)^{rp}\le\frac{1}{2^{rp}c^p}$.

This implies immediately:
$$\forall\;j\in\NN\;\forall\;t\in\RR:\;\;\sum_{p=0}^{\infty}\left|a_p\chi^{(rj)}_{cl,rp}(t)\right|\le 2|\mathbf{a}|_{M,c}c^jd^{rj}N_j.$$
So we can define the extension map $T_c:\Lambda_{M,c}\rightarrow\mathcal{D}_{r,N,cd^r}([-1,1])$ by $\mathbf{a}\mapsto\sum_{p=0}^{\infty}a_p\chi_{cl,rp}$, since
$$(\mathcal{B}^r\circ T_{c})(a)=\left(\sum_{p=0}^{\infty}a_p\chi^{(rj)}_{cl,rp}(0)\right)_{j\in\NN}=\left(\sum_{p=0}^{\infty}a_p\delta_{p,j}\right)_{j\in\NN}=(a_j)_{j\in\NN}=\mathbf{a}.$$
Note that the number $d$ is not depending on chosen $c$ and that $\sum_{p=0}^{\infty}a_p\chi^{(rj+k)}_{cl,rp}(0)=0$ for all $k\in\{1,\dots,r-1\}$.
\qed\enddemo

Now we are going to prove the second half of Theorem \ref{Roumieu-Surjectivitytheorem} and recall that by $(iii)$ in Remark \ref{Roumieu-Surjectivitytheoremremark} the inclusion $\Lambda_{\{M\}}\subseteq\{(f^{(rj)}(0))_{j\in\NN}: f\in\mathcal{D}_{r,\{N\}}([-1,1])\}$ does imply the nontriviality of $\mathcal{D}_{r,\{N\}}([-1,1])$ and condition \hyperlink{mnqr}{$(\text{nq}_r)$} for $N$.\vspace{6pt}

First we are proving the next result which generalizes \cite[Proposition 3.3]{surjectivity} and we follow the lines of the proof there.

\begin{proposition}\label{Roumieu-Surjectivitytheorem2}
Assume that
$$\Lambda_{\{M\}}\subseteq\{(f^{(rj)}(0))_{j\in\NN}: f\in\mathcal{D}_{r,\{N\}}([-1,1])\}$$
holds. Then for all $s\in\NN_{>0}$ there exists $h\in\NN_{>0}$ such that there exists a continuous linear extension map $T_s$ from $\Lambda_{M,s}$ into $\mathcal{D}_{r,N,h}([-1,1])$.
\end{proposition}

\demo{Proof}
For $j\in\NN$ consider the continuous linear functional $\tau^{j}: \mathcal{D}_{r,\{N\}}([-1,1])\rightarrow\RR$, $f\mapsto f^{(rj)}(0)$. Let
$$H:=\{f\in\mathcal{D}_{r,\{N\}}([-1,1]): f^{(rj)}(0)=0\;\forall\;j\in\NN\},$$
so $H$ is a closed subspace of $\mathcal{D}_{r,\{N\}}([-1,1])$. Let $\psi:\mathcal{D}_{r,\{N\}}([-1,1])\rightarrow\mathcal{D}_{r,\{N\}}([-1,1])/H$ be the canonical surjection.

Let $\mathbf{a}=(a_p)_p\in\Lambda_{\{M\}}$ be given, then there exists a function $f_{\mathbf{a}}\in\mathcal{D}_{r,\{N\}}([-1,1])$ such that $f_{\mathbf{a}}^{(rj)}(0)=a_j$ for all $j\in\NN$. We introduce the (well-defined) linear mapping $\phi:\Lambda_{\{M\}}\rightarrow\mathcal{D}_{r,\{N\}}([-1,1])/H$, $\mathbf{a}\mapsto\psi(f_{\mathbf{a}})$.

{\itshape Claim.} $\phi$ is continuous. Both $\Lambda_{\{M\}}$ and $\mathcal{D}_{r,\{N\}}([-1,1])/H$ are countable $(LB)$-spaces, so it suffices to prove that the graph of $\phi$ is sequentially closed. Consider a sequence $(\mathbf{a}^l)_{l\in\NN}$ in $\Lambda_{\{M\}}$ such that $\mathbf{a}^l\rightarrow\mathbf{a}$ as $l\rightarrow\infty$ and such that $(\phi(\mathbf{a}^l))_{l\in\NN}=(\psi(f_{\mathbf{a}^l}))_{l\in\NN}$ converges to $\psi(f)$ in $\mathcal{D}_{r,\{N\}}([-1,1])/H$.

Clearly for all $j\in\NN$ we have $a^l_j\rightarrow a_j$ as $j\rightarrow\infty$. Since each mapping $\tau^{j}$ is vanishing on $H$, the mapping $\widetilde{\tau^{j}}$ obtained by $\widetilde{\tau^{j}}\circ\psi=\tau^j$ is a continuous linear functional on $\mathcal{D}_{r,\{N\}}([-1,1])/H$ such that
$$\forall\;j\in\NN:\;a^l_j=\widetilde{\tau^{j}}(\phi(\mathbf{a}^l))\longrightarrow\widetilde{\tau^{j}}(\psi(f))=f^{(rj)}(0),$$
as $l\rightarrow\infty$. So $a_j=f^{(rj)}(0)$ for all $j\in\NN$, i.e. $\phi(\mathbf{a})=\psi(f)$ which proves the claim.\vspace{6pt}

Since we are dealing with two countable $(LB)$-spaces we can apply Grothendieck's factorization theorem (e.g. see \cite[24.33]{meisevogt}) to obtain
\begin{equation}\label{phipsi}
\forall\;s\in\NN_{>0}\;\exists\;h\in\NN_{>0}:\;\phi(\Lambda_{M,2s})\subseteq\psi(\mathcal{D}_{r,N,h}([-1,1]))=:E.
\end{equation}
We endow $E$ with the Banach space structure coming from its canonical identification with

$\mathcal{D}_{r,N,h}([-1,1])/(H\cap\mathcal{D}_{r,N,h}([-1,1]))$ and denote its norm by $\|\cdot\|_E$.\vspace{6pt}

So the mapping $\phi:\Lambda_{M,2s}\rightarrow E$ is continuous and linear between Banach spaces, hence
$$\exists\;C>0\;\forall\;\mathbf{a}\in\Lambda_{M,2s}:\;\;\|\phi(\mathbf{a})\|_E\le C|\mathbf{a}|_{M,2s}.$$
For $p\in\NN$ let $e^p:=(\delta_{p,j})_{j\in\NN}$, the $p$-th unit vector. We have $|e^p|_{M,2s}=\frac{1}{(2s)^p M_p}$ and so clearly $e^p\in\Lambda_{M,2s}$ for each $p\in\NN$. By \eqref{phipsi} there exists $\chi_p\in\mathcal{D}_{r,N,h}([-1,1])$ such that $\phi(e^p)=\psi(\chi_p)$ and $\|\chi_p\|_{r,N,h}\le2\|\phi(e^p)\|_E$. For this last inequality recall that $\|\cdot\|_E$ is the norm on the quotient space $E$ with $\|\phi(e^p)\|_E=\inf\{\|f\|_{r,N,h}: f\in\phi(e^p)=\psi(\chi_p)\}$.

Hence we get
\begin{align*}
|a_pe^p|_{M,2s}&\le|a_p||e^p|_{M,2s}\le|\mathbf{a}|_{M,s} s^p M_p\frac{1}{(2s)^p M_p}=2^{-p}|\mathbf{a}|_{M,s},
\end{align*}
and summarizing
\begin{align*}
\|a_p\chi_p\|_{r,N,h}&\le 2\|\phi(a_p e^p)\|_E\le2C|a_pe^p|_{M,2s}\le 2C2^{-p}|\mathbf{a}|_{M,s}.
\end{align*}
Consider now the series $T_s(\mathbf{a}):=\sum_{p=0}^{\infty}a_p\chi_p$. It defines a linear extension mapping $T_s:\Lambda_{M,s}\rightarrow\mathcal{D}_{r,N,h}([-1,1])$ and is continuous since
\begin{align*}
\|T_s(\mathbf{a})\|_{r,N,h}&\le\sum_{p=0}^{\infty}\|a_p\chi_p\|_{r,N,h}\le 2C|\mathbf{a}|_{M,s}\sum_{p=0}^{\infty}2^{-p}=4C|\mathbf{a}|_{M,s}
\end{align*}
for all $\mathbf{a}\in\Lambda_{M,s}$. Finally we have
$$(\mathcal{B}^r\circ T_s)(\mathbf{a})=\left(\sum_{p=0}^{\infty}a_p\chi_p^{(rj)}(0)\right)_{j\in\NN}=\left(\sum_{p=0}^{\infty}a_pe^p_j\right)_{j\in\NN}=\left(\sum_{p=0}^{\infty}a_p\delta_{p,j}\right)_{j\in\NN}=(a_j)_{j\in\NN}=\mathbf{a}.$$
\qed\enddemo

For the next theorem we need the following result, see \cite[Lemma 1.3.6]{hoermander}.

\begin{lemma}\label{Roumieu-Surjectivitytheorem3}
Let $l\in\NN_{>0}$ and $a_1,\dots,a_l$ be a nonincreasing sequence of positive real numbers with $0<A\le\sum_{i=1}^l a_i$. Then for all $f\in\mathcal{C}^l((-\infty,A])$ vanishing on $(-\infty,0]$ we have
$$\forall\;t\le A:\hspace{30pt}|f(t)|\le\sum_{j\in J_l}2^{2j} a_1\cdots a_j\sup_{s<t}|f^{(j)}(s)|,$$
where $J_l:=\{j\in\NN_{>0}: 1\le j\le l, a_{j+1}<a_j\;\text{or}\;j=l\}$.
\end{lemma}

Our next result proves the converse statement of Theorem \ref{Roumieu-Surjectivitytheorem1} and generalizes \cite[Theorem 3.5]{surjectivity}.

\begin{theorem}\label{Roumieu-Surjectivitytheorem4}
If $\Lambda_{\{M\}}\subseteq\{(f^{(rj)}(0))_{j\in\NN}: f\in\mathcal{D}_{r,\{N\}}([-1,1])\}$ is valid, then \hyperlink{SV}{$(M,N)_{SV_r}$} holds true.
\end{theorem}

\demo{Proof}
By Proposition \ref{Roumieu-Surjectivitytheorem2} there exist $s\in\NN_{>0}$ and a continuous linear extension map $T$ from $\Lambda_{M,1}$ into $\mathcal{D}_{r,N,s}([-1,1])$. We choose $D>\max\{1,\|T\|\}$ and $1>h>0$ to be small enough to guarantee $0<4hs^2<1/2$. For $p\in\NN_{>0}$ we consider the (increasing) sequence $\tau^p:=(\tau^p_j)_{j\ge 1}$ defined by
$$\tau^p=\big(\underbrace{\left(\frac{\nu_{2p}}{h}\right)^{1/r},\dots,\left(\frac{\nu_{2p}}{h}\right)^{1/r}}_{rp-\text{times}},\underbrace{\left(\frac{\nu_{2p+1}}{h}\right)^{1/r},\dots,\left(\frac{\nu_{2p+1}}{h}\right)^{1/r}}_{r-\text{times}},\underbrace{\left(\frac{\nu_{2p+2}}{h}\right)^{1/r},\dots,\left(\frac{\nu_{2p+2}}{h}\right)^{1/r}}_{r-\text{times}},\dots\big).$$
Put $e^p:=(\delta_{p,j})_{j\in\NN}$ and so $e^p\in\Lambda_{M,1}$ with $|e^p|_{M,1}=\frac{1}{M_p}$ for each $p\in\NN_{>0}$. Moreover put $\chi_p:=T(e^p)$, so $\chi_p\in\mathcal{D}_{r,N,s}([-1,1])$ and $\chi^{(rj)}_p(0)=e^p_j=\delta_{p,j}$. Finally, for all $p\in\NN_{>0}$ and $0\le j\le p-1$ we introduce the function
\begin{equation*}\label{varrhoroumieu}
\varrho_{p,j}(t):=
\begin{cases}
0 &\text{if}\;t\le 0,
\\
\chi_p^{(rj)}(t)-\frac{t^{r(p-j)}}{(r(p-j))!} &\text{if}\;t>0.
\end{cases}
\end{equation*}
We want to apply Lemma \ref{Roumieu-Surjectivitytheorem3} and first note that $\varrho_{p,j}$ is smooth on $\RR$. Of course it is smooth on $(-\infty,0)$ and on $(0,+\infty)$ and since for all $k\in\NN$, $\chi^{(rj+k)}_p(0)=\delta_{p,rj+k}=1$ if and only if $r(p-j)=k$, it is also smooth at $0$. Let $z$ with
\begin{equation}\label{z}
0<z<\frac{rph}{(\nu_{2p})^{1/r}}+r\sum_{l=2p+1}^{\infty}\left(\frac{h}{\nu_l}\right)^{1/r},
\end{equation}
and use Lemma \ref{Roumieu-Surjectivitytheorem3} for the function $\varrho_{p,j}$ and the sequence $\tau^p$, i.e. put $a_j:=(\tau^p_j)^{-1}$. Hence we get for each $z$ as in \eqref{z}:
\begin{equation}\label{Roumieutheorem4equ1}
|\varrho_{p,j}(z)|\le\sum_{k=p}^{\infty}\frac{(4h)^k}{(\nu_{2p})^{p}\prod_{j={p+1}}^k\nu_{p+j}}\sup_{t\in[0,z]}|\varrho^{(rk)}_{p,j}(t)|.
\end{equation}
Concerning the index in the summation we recall that we can start at $k=rp$, since $\tau^p_j=\tau^p_{j+1}$ for $1\le j\le rp-1$, and moreover $\tau^p_j=\tau^p_{j+1}$ for $r(p+i)<j<j+1\le r(p+i+1)$ for any $i\in\NN$. Then we estimate as follows for $k\ge p\Leftrightarrow rk\ge rp$:
\begin{align*}
\sup_{t\in[0,z]}|\varrho^{(rk)}_{p,j}(t)|&\le\sup_{t\in[0,z]}|\chi^{(r(j+k))}_{p}(t)|+1\le\|\chi_p\|_{r,N,s}s^{j+k}N_{j+k}+1
\\&
\underbrace{\le}_{(\star)}D|e^p|_{M,1}s^{2k}N_{j+k}+1\le D s^{2k}\frac{N_{j+k}}{M_p}+1\le 2Ds^{2k}\frac{N_{j+k}}{M_p}.
\end{align*}
$(\star)$ holds since $\|\chi_p\|_{r,N,s}\le\|T\||e^p|_{M,1}\le D|e^p|_{M,1}$ and since $s^{j+k}\le s^{2k}$ by $0\le j<p\le k$. Moreover $p\le k$ implies $M_p\le N_p\le N_k\le N_{j+k}$ which was used for the last estimate.

On the other hand we have for all $p\in\NN_{>0}$, $0\le j\le p-1$ and $k\ge p$:
\begin{equation}\label{Roumieutheorem4equ2}
\frac{1}{(\nu_{2p})^p\prod_{j={p+1}}^k\nu_{p+j}}\le\frac{N_j}{N_{j+k}},
\end{equation}
which holds because
\begin{align*}
N_{j+k}&=N_j\nu_{j+1}\cdots\nu_{j+k}\le N_j\nu_{p+1}\cdots\nu_{p+k}\le N_j(\nu_{2p})^p\prod_{j={p+1}}^k\nu_{p+j}.
\end{align*}
Using these estimates we are going to prove now for each $z$ as in \eqref{z}:
\begin{equation}\label{Roumieutheorem4equ3}
|\varrho_{p,j}(z)|\le 2D\frac{N_j}{M_p}2^{-p+1}.
\end{equation}
This holds by the following calculation and the sufficient small choice of $h$:
\begin{align*}
|\varrho_{p,j}(z)|&\underbrace{\le}_{\eqref{Roumieutheorem4equ1},\eqref{Roumieutheorem4equ2}}\sum_{k=p}^{\infty}(4h)^k\frac{N_j}{N_{j+k}}\sup_{t\in[0,z]}|\varrho^{(rk)}_{p,j}(t)|\le\sum_{k=p}^{\infty}(4h)^k\frac{N_j}{N_{j+k}}2Ds^{2k}\frac{N_{j+k}}{M_p}
\\&
\le 2D\frac{N_j}{M_p}\sum_{k=p}^{\infty}(\underbrace{4hs^2}_{<1/2})^k\le 2D\frac{N_j}{M_p}2^{-p+1}.
\end{align*}
In the next step we prove for any $z>0$ as in \eqref{z} that
\begin{equation}\label{Roumieutheorem4equ4}
\forall\;p\in\NN_{>0}\;\forall\;j\in\NN, 0\le j\le p-1:\;\;\frac{1}{2}\frac{z^{r(p-j)}}{(r(p-j))!}\le 2D\frac{s^pN_j}{M_p}.
\end{equation}
Let $p\in\NN_{>0}$ and $0\le j\le p-1$, we distinguish two cases:

{\itshape Case $1$.} If $\chi_p^{(rj)}(z)\le\frac{1}{2}\frac{z^{r(p-j)}}{(r(p-j))!}$, then
\begin{align*}
\frac{1}{2}\frac{z^{r(p-j)}}{(r(p-j))!}\le\frac{z^{r(p-j)}}{(r(p-j))!}-\chi_p^{(rj)}(z)=|\varrho_{p,j}(z)|\le 2D\frac{s^p N_j}{M_p}.
\end{align*}
For the last step we have used \eqref{Roumieutheorem4equ3} and $2^{-p+1}\le s^p\Leftrightarrow 2\le(2s)^p$.

{\itshape Case $2$.} If $\chi_p^{(rj)}(z)>\frac{1}{2}\frac{z^{r(p-j)}}{(r(p-j))!}$, then
\begin{align*}
\frac{1}{2}\frac{z^{r(p-j)}}{(r(p-j))!}&<\chi_p^{(rj)}(z)\le\|\chi_p\|_{r,N,s}s^jN_j\le\|T\||e^p|_{M,1}s^j N_j\le D\frac{s^p N_j}{M_p}.
\end{align*}
For the last inequality we have also used $s^j\le s^{p}$ since $j\le p$. Hence, by choosing $z=\sum_{l=2p+1}^{\infty}\left(\frac{h}{\nu_l}\right)^{1/r}$ which is possible by \eqref{z} and having \hyperlink{mnqr}{$(\text{nq}_r)$} for $N$, we obtain for all $p\in\NN_{>0}$ and $0\le j\le p-1$
\begin{equation*}\label{importantsqequ}
\sum_{l=2p+1}^{\infty}\left(\frac{h}{\nu_l}\right)^{1/r}=z\underbrace{\le}_{\eqref{Roumieutheorem4equ4}}(4D)^{1/(r(p-j))}(r(p-j))!^{1/(r(p-j))}\left(\frac{s^p N_j}{M_p}\right)^{1/(r(p-j))}\le 4Drp\left(\frac{s^p N_j}{M_p}\right)^{1/(r(p-j))},
\end{equation*}
where the last inequality holds since $(r(p-j))!^{1/(r(p-j))}\le rp\Leftrightarrow(r(p-j))!\le(rp)^{r(p-j)}$ for all $0\le j\le p-1$. By definition of the expression $\lambda^{M,N}_{p,s}$ in \eqref{lambdaexpressionequ} (write again simply $\lambda_{p,s}$ for it) we have shown so far
\begin{equation}\label{Roumieutheorem4equ5}
\forall\;p\in\NN_{>0}:\;\;\frac{(\lambda_{p,s})^{1/r}}{p}\sum_{l=2p+1}^{\infty}\left(\frac{1}{\nu_l}\right)^{1/r}\le\frac{4Dr}{h^{1/r}}.
\end{equation}
Now we are able to show \hyperlink{SV}{$(M,N)_{SV_r}$} and finish the proof. Let $p\in\NN_{>0}$ be arbitrary, then
\begin{align*}
\frac{(\lambda_{p,s})^{1/r}}{p}\sum_{l=p}^{\infty}\left(\frac{1}{\nu_l}\right)^{1/r}\underbrace{\le}_{\eqref{Roumieutheorem4equ5}}\frac{(\lambda_{p,s})^{1/r}}{p}\sum_{l=p}^{2p}\left(\frac{1}{\nu_l}\right)^{1/r}+\frac{4Dr}{h^{1/r}}\underbrace{\le}_{\eqref{lambdaexpressionequ1}}\frac{(\nu_p)^{1/r}}{p}\frac{p+1}{(\nu_p)^{1/r}}+\frac{4Dr}{h^{1/r}}\le 2+\frac{4Dr}{h^{1/r}}.
\end{align*}
\qed\enddemo

\section{The image of the Borel mapping in the Beurling case}\label{section5}
Let from now on in this section $M,N\in\RR_{>0}^{\NN}$ and $r\in\NN_{\ge 1}$ be such that $(I)$ and $(III)$ from Section \ref{section4} are valid and moreover

\begin{itemize}
\item[$(II)_{B,r}$] $\lim_{p\rightarrow\infty}\left(\frac{(M_p)^{1/r}}{p!}\right)^{1/p}=\infty$ (the letter $B$ in the notation stands for {\itshape Beurling}).
\end{itemize}

Analogously as in Remark \ref{roumieucondtionremark} in the Roumieu case above we have:

\begin{remark}\label{beurlingcondtionremark}
\begin{itemize}
\item[$(i)$] Again we can assume without loss of generality that even $M\le N$ holds true.

\item[$(ii)$] For any $r>1$ condition $(II)_{B,r}$ implies $(II)_{B,1}$, i.e. $\lim_{p\rightarrow\infty}(m_p)^{1/p}=\infty$, because $\frac{(M_p)^{1/r}}{p!}=m_p(M_p)^{\frac{1}{r}-1}$.

\item[$(iii)$] $(II)_{B,r}$ for $M$ does imply $(II)_{B,1}$ for $P^{M,r}$ (here for all $C\ge 1$ large there exists $p_C$ such that for all $p\ge p_C$ we have $M_p\ge C^{rp}p!^r$).
\end{itemize}
\end{remark}

The goal of this section is to prove the following characterization.

\begin{theorem}\label{Beurling-Surjectivitytheorem}
Let $M$ and $N$ be as assumed above and $r\in\NN_{\ge 1}$. Then
\begin{equation}\label{Beurling-Surjectivitytheoremequ}
\Lambda_{(M)}\subseteq\{(f^{(rj)}(0))_{j\in\NN}: f\in\mathcal{D}_{r,(N)}([-1,1])\}
\end{equation}
holds if and only if \hyperlink{SV}{$(M,N)_{SV_r}$} is satisfied.
\end{theorem}

\begin{remark}\label{Beurling-Surjectivitytheoremremark}
Similarly as in the Roumieu case above, this result is extending \cite[Theorem 1.1]{surjectivity} also for general $r$-interpolating spaces, and even in the case $r=1$ our approach is slightly stronger than the result from \cite{surjectivity} since

\begin{itemize}
\item[$(i)$] we only require $M\hyperlink{mpreceq}{\preceq}N$ instead of the stronger assumption $\mu_p\le\nu_p$ for all $p\in\NN$ and

\item[$(ii)$] assumption \hyperlink{mnqr}{$(\on{nq}_r)$} for $N$ is not needed because even in the general setting analogously as commented in Remark \ref{Roumieu-Surjectivitytheoremremark} above the inclusion \eqref{Beurling-Surjectivitytheoremequ} does imply that $\mathcal{D}_{r,(N)}([-1,1])$ is nontrivial and Theorem \ref{nonquasitheorem} yields \hyperlink{mnqr}{$(\on{nq}_r)$} for $N$.
\end{itemize}
\end{remark}

The strategy is to reduce the proof to the Roumieu case, as it has been done in \cite[Section 4]{surjectivity}, and to do so we will have to apply the following result, see \cite[Lemme 16]{ChaumatChollet94}.

\begin{lemma}\label{Beurling-Surjectivitytheorem1}
Let $(\alpha_k)_{k\in\NN_{>0}}$ be a sequence of non-negative real numbers such that $\sum_{k=1}^{\infty}\alpha_k<\infty$. Furthermore let $\beta=(\beta_k)_{k\in\NN_{>0}}$ and $\gamma=(\gamma_k)_{k\in\NN_{>0}}$ be sequences of positive real numbers such that $\lim_{k\rightarrow\infty}\beta_k=0=\lim_{k\rightarrow\infty}\gamma_k$, and assume $\gamma$ is nonincreasing.

Then there exists a sequence $(\lambda_k)_{k\in\NN_{>0}}$ such that
\begin{itemize}
\item[$(i)$] $k\mapsto\lambda_k$ is nondecreasing,
\item[$(ii)$] $\lim_{k\rightarrow\infty}\lambda_k=\infty$,
\item[$(iii)$] $k\mapsto\lambda_k\gamma_k$ is nonincreasing,
\item[$(iv)$] $\lim_{k\rightarrow\infty}\lambda_k\beta_k=0$,
\item[$(v)$] $\forall\;p\in\NN_{>0}:\;\;\sum_{k=p}^{\infty}\lambda_k\alpha_k\le 8\lambda_p\sum_{k=p}^{\infty}\alpha_k$.
\end{itemize}
\end{lemma}

The next result generalizes \cite[Theorem 4.2.]{surjectivity}.

\begin{theorem}\label{Beurling-Surjectivitytheorem2}
Let $M$ and $N$ be as assumed above, $r\in\NN_{\ge 1}$ and satisfying \hyperlink{SV}{$(M,N)_{SV_r}$}. Then we have
$$\Lambda_{(M)}\subseteq\{(f^{(rj)}(0))_{j\in\NN}: f\in\mathcal{D}_{r,(N)}([-1,1])\}.$$
\end{theorem}

\demo{Proof}
Let $0\neq\mathbf{a}=(a_p)_p\in\Lambda_{(M)}$, so by definition for all $h\in\NN_{>0}$ there exists $C_h>0$ such that for all $p\in\NN$ we have $|a_p|\le C_h h^{-p} M_p$, hence $\left(\frac{|a_p|}{M_p}\right)^{1/p}\le\frac{C_h^{1/p}}{h}$ for all $p\in\NN$ and $h\in\NN_{>0}$ which gives
\begin{equation}\label{Beurling-equ1}
\lim_{p\rightarrow\infty}\left(\frac{|a_p|}{M_p}\right)^{1/p}=0.
\end{equation}
Now define a sequence $(\epsilon_p)_{p\in\NN_{>0}}$ by $\epsilon_p:=\sup_{k\ge p}\left(\frac{|a_k|}{M_k}\right)^{1/k}$, which is clearly nonincreasing and $\lim_{p\rightarrow\infty}\epsilon_p=0$. Since by definition $\epsilon_j\ge\left(\frac{|a_p|}{M_p}\right)^{1/p}$ whenever $1\le j\le p$ we also get $|a_p|\le\epsilon_1\cdots\epsilon_p M_p$ for all $p\in\NN_{>0}$.

Put $\alpha_k:=0$, $\beta_k:=\max\{\epsilon_k,\frac{k}{(M_k)^{1/(rk)}}\}$ and $\gamma_k:=\frac{1}{\mu_k}$ for all $k\in\NN_{>0}$. By \eqref{Beurling-equ1}, standard assumption $(II)_{B,r}$ on $M$ and Stirling's formula we have $\beta_k\rightarrow 0$ as $k\rightarrow\infty$. $\gamma=(\gamma_k)_k$ is nonincreasing and tending to $0$ by \eqref{mucompare}. So we can apply Lemma \ref{Beurling-Surjectivitytheorem1} to obtain a sequence $\theta=(\theta_k)_{k\in\NN_{>0}}$ such that $\theta$ is nondecreasing, tending to $\infty$, $\theta_k\gamma_k$ is nonincreasing and finally $\theta_k\beta_k\rightarrow 0$ as $k\rightarrow\infty$. W.l.o.g. we can assume $\theta_1=1$ and moreover we have $\theta_k\gamma_k\rightarrow 0$, since
$$\forall\;k\in\NN_{>0}:\;\;\theta_k\gamma_k=\frac{\theta_k}{\mu_k}\underbrace{\le}_{(M_k)^{1/r}\le M_k\le(\mu_k)^k}\frac{\theta_k}{(M_k)^{1/(rk)}}\le\theta_k\frac{k}{(M_k)^{1/(rk)}}\underbrace{\le}_{\text{def. of}\;\beta}\theta_k\beta_k.$$
In the next step we apply Lemma \ref{Beurling-Surjectivitytheorem1} to $\alpha'_k=\gamma'_k:=\left(\frac{1}{\nu_k}\right)^{1/r}$ and $\beta'_k:=\max\{\frac{1}{(\theta_{\lfloor k/2\rfloor})^{1/(2r)}},\frac{1}{(\nu_k)^{1/r}}\}$ for $k\in\NN_{>0}$, where $\lfloor k/2\rfloor$ denotes the integer part of $k/2$ and we put $\theta_0:=1$. This can be done since $(\nu_k)^{1/r}\rightarrow\infty$ as $k\rightarrow\infty$ by having $\sum_{k\ge 1}\left(\frac{1}{\nu_k}\right)^{1/r}<\infty$ (\hyperlink{mnqr}{$(\text{nq}_r)$} for $N$ follows by assumption \hyperlink{SV}{$(M,N)_{SV_r}$}) and since $\theta_k\rightarrow\infty$ as $k\rightarrow\infty$.

Hence we obtain another sequence $\theta'=(\theta'_p)_{p\in\NN_{>0}}$ which is nondecreasing, tending to infinity, $\theta'_k\gamma'_k$ is nonincreasing, $\theta'_k\beta'_k\rightarrow 0$ and finally
\begin{equation}\label{Beurling-equ2}
\forall\;p\in\NN_{>0}:\;\;\;\sum_{k=p}^{\infty}\frac{\theta'_k}{(\nu_k)^{1/r}}\le 8\theta'_p\sum_{k=p}^{\infty}\left(\frac{1}{\nu_k}\right)^{1/r}.
\end{equation}
W.l.o.g. we can assume $\theta'_1=1$. By definition and since $\theta'_k\beta'_k\rightarrow 0$ we get $\lim_{k\rightarrow\infty}\frac{\theta'_k}{(\theta_{\lfloor k/2\rfloor})^{1/(2r)}}=0=\lim_{k\rightarrow\infty}\frac{\theta'_k}{(\nu_k)^{1/r}}$, which shows $\theta'_k\gamma'_k\rightarrow 0$ as $k\rightarrow\infty$. Since $\theta$ is nondecreasing, $\theta_p\ge 1$ for each $p$ we also obtain
\begin{equation}\label{Beurling-equ3}
\exists\;A\ge 1\;\forall\;k\in\NN_{>0}:\;\;\;(\theta'_k)^r\le A(\theta_{\lfloor k/2\rfloor})^{1/2}\le A\theta_{\lfloor k/2\rfloor}\le A\theta_k.
\end{equation}
Using \eqref{Beurling-equ3} we are going to show:
\begin{equation}\label{Beurling-equ4}
\exists\;A\ge 1\;\forall\;p\in\NN_{>0}\;\forall\;j\in\NN,0\le j\le p-1:\;\;\;\frac{(\theta'_p)^r}{(\theta_{j+1}\cdots\theta_p)^{1/(p-j)}}\le A.
\end{equation}
For $0\le j\le\lfloor p/2\rfloor-1$ we have
$$(\theta_{j+1}\cdots\theta_p)^{1/(p-j)}\ge(\theta_{\lfloor p/2\rfloor}\cdots\theta_p)^{1/(p-j)}\ge(\theta_{\lfloor p/2\rfloor})^{(p-\lfloor p/2\rfloor)/(p-j)}\ge(\theta_{\lfloor p/2\rfloor})^{(p-\lfloor p/2\rfloor)/p}\ge\sqrt{\theta_{\lfloor p/2\rfloor}},$$
since $\theta$ is nondecreasing, $\theta_p\ge 1$ for each $p$ and for the last inequality we have used $\frac{p-\lfloor p/2\rfloor}{p}=1-\frac{\lfloor p/2\rfloor}{p}\ge\frac{1}{2}\Leftrightarrow\frac{p}{2}\ge\lfloor p/2\rfloor$. On the other hand, if $\lfloor p/2\rfloor\le j\le p-1$, then $(\theta_{j+1}\cdots\theta_p)^{1/(p-j)}\ge\theta_{j+1}\ge\theta_{\lfloor p/2\rfloor}$. Summarizing we end up with $\frac{(\theta'_p)^r}{(\theta_{j+1}\cdots\theta_p)^{1/(p-j)}}\le\frac{(\theta'_p)^r}{\sqrt{\theta_{\lfloor p/2\rfloor}}}$ respectively $\frac{(\theta'_p)^r}{(\theta_{j+1}\cdots\theta_p)^{1/(p-j)}}\le\frac{(\theta'_p)^r}{\theta_{\lfloor p/2\rfloor}}$, now apply \eqref{Beurling-equ3}.\vspace{6pt}

We introduce sequences $R=(R_p)_p$ and $S=(S_p)_p$ defined by $R_k:=\prod_{j=0}^k\varrho_j$, $S_k:=\prod_{j=0}^k\sigma_j$, where $$\forall\;k\in\NN_{>0}:\;\;\varrho_k:=\frac{\mu_k}{\theta_k},\hspace{20pt}\sigma_k:=A\frac{\nu_k}{(\theta'_k)^r},\hspace{60pt}\varrho_0=\sigma_0:=1.$$
We summarize:

\begin{itemize}
\item[$(i)$] $R$ is log-convex since $(\varrho_k)^{-1}=\frac{\theta_k}{\mu_k}=\theta_k\gamma_k$ is nonincreasing (and tending to $0$ as $k\rightarrow\infty$).
\item[$(ii)$] $S$ is log-convex since $(\sigma_k)^{-1}=A^{-1}\frac{(\theta'_k)^r}{\nu_k}=A^{-1}(\theta'_k\gamma'_k)^r$ is nonincreasing (and tending to $0$ as $k\rightarrow\infty$).
\item[$(iii)$] $R_k\le S_k$ since $R_k=M_k(\theta_1\cdots\theta_k)^{-1}\le N_k A^k(\theta'_1\cdots\theta'_k)^{-r}=S_k$ for each $k\in\NN_{>0}$, see \eqref{Beurling-equ3}.
\item[$(iv)$] $S$ satisfies \hyperlink{mnqr}{$(\text{nq}_r)$} since by definition, \hyperlink{mnqr}{$(\text{nq}_r)$} for $N$ and \eqref{Beurling-equ2} we get
    \begin{equation}\label{Beurling-equ5}
    \forall\;p\in\NN_{>0}:\;\;\;\sum_{k\ge p}\left(\frac{1}{\sigma_k}\right)^{1/r}=A^{-1/r}\sum_{k\ge p}\frac{\theta'_k}{(\nu_k)^{1/r}}\le A^{-1/r}8\theta'_p\sum_{k\ge p}\left(\frac{1}{\nu_k}\right)^{1/r}<\infty.
    \end{equation}
\item[$(v)$] $R$ satisfies $\lim_{k\rightarrow\infty}\frac{(R_k)^{1/(rk)}}{(k!)^{1/k}}=\infty$, i.e. $(II)_{B,r}$: For all $k\in\NN_{>0}$ we get
    $$\frac{k}{(R_k)^{1/(rk)}}=\frac{k}{((\theta_1\cdots\theta_k)^{-1}M_k)^{1/(rk)}}\le\frac{k(\theta_k)^{1/r}}{(M_k)^{1/(rk)}}\le\frac{k\theta_k}{(M_k)^{1/(rk)}}\le\theta_k\beta_k,$$
    which tends to $0$ as $k\rightarrow\infty$. We have used that $\theta$ is nondecreasing with $\theta_1=1$, and finally $\beta_k\ge\frac{k}{(M_k)^{1/k}}$ by definition of $\beta$. The conclusion follows by applying Stirling's formula.
\item[$(vi)$] $\lim_{k\rightarrow\infty}(S_k)^{1/k}=\infty$ holds by $(iii)$ and $(v)$, in fact we even have $\lim_{k\rightarrow\infty}\frac{(S_k)^{1/(rk)}}{(k!)^{1/k}}=\infty$.
\end{itemize}
\vspace{6pt}

Let $p\in\NN_{>0}$ and $0\le j\le p-1$ be arbitrary, then we calculate as follows by using \hyperlink{SV}{$(M,N)_{SV_r}$} in the second estimate below:
\begin{align*}
\sum_{k=p}^{\infty}\left(\frac{1}{\sigma_k}\right)^{1/r}&\underbrace{\le}_{\eqref{Beurling-equ5}}A^{-1/r}8\theta'_p\sum_{k=p}^{\infty}\left(\frac{1}{\nu_k}\right)^{1/r}\le A^{-1/r}8p\theta'_pD\left(\frac{s^p\nu_1\cdots\nu_j}{\mu_1\cdots\mu_p}\right)^{1/(r(p-j))}
\\&
=A^{-1/r}8p\theta'_pD\left(\frac{s^p(\theta'_1\cdots\theta'_j)^r\sigma_1\cdots\sigma_j}{A^j\theta_1\cdots\theta_p\varrho_1\cdots\varrho_p}\right)^{1/(r((p-j))}=A^{-1/r}8p\theta'_pD\left(\frac{s^p(\theta'_1\cdots\theta'_j)^rS_j}{A^j\theta_1\cdots\theta_pR_p}\right)^{1/(r(p-j))}
\\&
\underbrace{\le}_{\eqref{Beurling-equ3}} A^{-1/r}8p\theta'_pD\left(\frac{s^pS_j}{\theta_{j+1}\cdots\theta_pR_p}\right)^{1/(r(p-j))}\underbrace{\le}_{\eqref{Beurling-equ4}}8Dp\left(\frac{s^pS_j}{R_p}\right)^{1/(r(p-j))}.
\end{align*}
So we have shown that
\begin{equation*}\label{Beurling-equ7}
\sup_{p\in\NN_{>0}}\frac{(\lambda^{R,S}_{p,s})^{1/r}}{p}\sum_{k=p}^{\infty}\left(\frac{1}{\sigma_k}\right)^{1/r}<\infty,
\end{equation*}
i.e. \hyperlink{SV}{$(R,S)_{SV_r}$}. $(i)-(vi)$ guarantee that all standard assumptions $(I)$, $(II)_{R,r}$ and $(III)$ in Section \ref{section4} on $R$ and $S$ are satisfied. By using Theorem \ref{Roumieu-Surjectivitytheorem1} we get
\begin{equation}\label{Beurling-equ8}
\Lambda_{\{R\}}\subseteq\{(f^{(rj)}(0))_{j\in\NN}: f\in\mathcal{D}_{r,\{S\}}([-1,1]).
\end{equation}
{\itshape Claim.} If $0\neq\mathbf{a}=(a_p)_p\in\Lambda_{(M)}$, then $\mathbf{a}\in\Lambda_{\{R\}}$. By definition for all $p\in\NN_{>0}$ we have
$$|a_p|\le\epsilon_1\cdots\epsilon_pM_p=\epsilon_1\theta_1\cdots\epsilon_p\theta_p(\theta_1\cdots\theta_p)^{-1}M_p=\epsilon_1\theta_1\cdots\epsilon_p\theta_p R_p.$$
Moreover $\epsilon_p\theta_p\le\beta_p\theta_p$ for all $p\in\NN_{>0}$ holds by definition of $\beta$, hence $\epsilon_p\theta_p\rightarrow 0$ as $p\rightarrow\infty$. So there exists $p_1\in\NN$ such that $|a_p|\le R_p$ for all $p\ge p_1$ which proves the claim.

Hence, given $\mathbf{a}\in\Lambda_{(M)}$, by applying Theorem \ref{Roumieu-Surjectivitytheorem1} there exists some $f_{\mathbf{a}}\in\mathcal{D}_{r,\{S\}}([-1,1])$ such that $f_{\mathbf{a}}^{(rj)}(0)=a_j$ for all $j\in\NN$ (see \eqref{Beurling-equ8}).\vspace{6pt}

{\itshape Claim.} $f_{\mathbf{a}}\in\mathcal{D}_{r,(N)}([-1,1])$. By assumption $\supp(f_{\mathbf{a}})\subseteq[-1,1]$ and there exists some $D>0$ and $h>0$ (both large) such that for all $n\in\NN$ and $x\in[-1,1]$ we have $|f_{\mathbf{a}}^{(rn)}(x)|\le D h^n S_n$ and finally $f_{\mathbf{a}}^{(rn+j)}(0)=0$ for all $j\in\{1,\dots,r-1\}$ and $n\in\NN$. Let $h_1\in\NN_{>0}$ be given (arbitrary large) but from now on fixed. Since $\theta'$ is nondecreasing and tending to infinity we can find $p_2\in\NN$ such that $\theta'_p\ge(h_1hA)^{1/r}$ for all $p\ge p_2$. So for all $p>p_2$ the following estimate is valid:
\begin{align*}
N_p&=\nu_1\cdots\nu_p=A^{-p}(\theta'_1\cdots\theta'_p)^r\sigma_1\cdots\sigma_p\ge A^{-p}(\theta'_{p_2+1}\cdots\theta'_p)^rS_p\ge A^{-p}((h_1hA)^{1/r})^{r(p-p_2)}S_p
\\&
=A^{-p_2}(hh_1)^{p-p_2}S_p,
\end{align*}
where in the first inequality we have used $\theta'_p\ge 1$ for each $p\in\NN$. Hence for $p>p_2$ we get:
$$|f_{\mathbf{a}}^{(rp)}(x)|\le D h^p S_p\le D h^p A^{p_2}(hh_1)^{p_2-p}N_p=D(Ahh_1)^{p_2}h_1^{-p}N_p,$$
which proves the second claim and finishes the proof.
\qed\enddemo

To prove the converse direction we use the notation introduced on \cite[p. 396]{surjectivity}. Let $s\in\NN_{>0}$, then we denote by $E_{s}$ the normed space $(\mathcal{D}_{r,(N)}([-1,1]),\|\cdot\|_{r,N,1/s})$ and by $F_{s}$ its completion. For $p\in\NN$ we consider the functional $\tau^{p}$ on $E_{s}$ defined by $\tau^{p}(f):=f^{(rp)}(0)$. It is continuous and linear and has a unique continuous linear extension on $F_{s}$ which will be still denoted by $\tau^{p}$.

$T_{s}:\Lambda_{(M)}\rightarrow F_{s}$ is called an {\itshape extension mapping} if $\tau^{p}(T(\mathbf{a}))=a_p$ for all $\mathbf{a}=(a_p)_p\in\Lambda_{(M)}$ and $p\in\NN$.

The next result generalizes \cite[Proposition 4.3]{surjectivity}.

\begin{proposition}\label{Beurling-Surjectivitytheorem3}
The inclusion
$$\Lambda_{(M)}\subseteq\{(f^{(rj)}(0))_{j\in\NN}: f\in\mathcal{D}_{r,(N)}([-1,1])\}$$
implies that for all $s\in\NN_{>0}$ (large) there exists a continuous linear extension mapping $T_{s}:\Lambda_{(M)}\rightarrow F_{s}$ such that $T_{s}(e^p)\in\mathcal{D}_{r,(N)}([-1,1])$ for all $p\in\NN$.
\end{proposition}

\demo{Proof}
Let $s\in\NN_{>0}$ be arbitrary (large) but from now on fixed and introduce
$$H:=\{f\in\mathcal{D}_{r,(N)}([-1,1]): f^{(rj)}(0)=0\;\forall\;j\in\NN\},$$
which is a closed sub-space of $\mathcal{D}_{r,(N)}([-1,1])$ and of $E_{s}$. Let $\psi:\mathcal{D}_{r,(N)}([-1,1])\rightarrow\mathcal{D}_{r,(N)}([-1,1])/H$ be the canonical surjection. Let $\mathbf{a}=(a_p)_p\in\Lambda_{(M)}$, then there exists some $f_{\mathbf{a}}\in\mathcal{D}_{r,(N)}([-1,1])$ such that $f_{\mathbf{a}}^{(rp)}(0)=a_p$ for all $p\in\NN$. As in Proposition \ref{Roumieu-Surjectivitytheorem2} we introduce the linear mapping $\phi:\Lambda_{(M)}\rightarrow\mathcal{D}_{r,(N)}([-1,1])/H$ defined by $\mathbf{a}\mapsto f_{\mathbf{a}}$, i.e. $\phi(\mathbf{a})=\psi(f_{\mathbf{a}})$. The continuity of $\phi$ follows analogously as in Proposition \ref{Roumieu-Surjectivitytheorem2} since by assumption $\phi$ is a linear mapping between two Fréchet spaces.

So it follows that $\phi$ is also linear and continuous between $\Lambda_{(M)}$ and $E_{s}/H$, i.e.
\begin{equation}\label{Beurling-Surjectivitytheorem3equ1}
\exists\;c\in\NN_{>0}:\;\exists\;A\ge 1\;\forall\;\mathbf{a}=(a_p)_p\in\Lambda_{(M)}:\;\;\|\phi(\mathbf{a})\|^{\widetilde{}}_{r,N,1/s}\le A|\mathbf{a}|_{M,1/c},
\end{equation}
where $\|\cdot\|^{\widetilde{}}$ shall denote the norm on the quotient $E_{s}/H$.

Let $p\in\NN$, then $e^p:=(\delta_{p,j})_{j\in\NN}$ belongs to $\Lambda_{(M)}$ with $|e^p|_{M,1/c}=\frac{c^p}{M_p}$ (for any $c\in\NN_{>0}$). Let $\chi_p\in E_{s}$ such that $\psi(\chi_p)=\phi(e^p)$ holds true and $\|\chi_p\|_{r,N,1/s}\le 2\|\phi(e^p)\|^{\widetilde{}}_{r,N,1/s}$. For $\mathbf{a}=(a_p)_p\in\Lambda_{(M)}$ we estimate as follows where $c$ is coming from \eqref{Beurling-Surjectivitytheorem3equ1}:
\begin{align*}
\left\|\sum_{p=0}^{\infty}a_p\chi_p\right\|^{\widehat{}}_{r,N,1/s}&\le\sum_{p=0}^{\infty}|a_p|\|\chi_p\|_{r,N,1/s}\le 2\sum_{p=0}^{\infty}|\mathbf{a}|_{M,1/(2c)}\frac{M_p}{(2c)^p}\|\phi(e^p)\|^{\widetilde{}}_{r,N,1/s}
\\&
\underbrace{\le}_{\eqref{Beurling-Surjectivitytheorem3equ1}}2A|\mathbf{a}|_{M,1/(2c)}\sum_{p=0}^{\infty}\frac{M_p}{(2c)^p}|e^p|_{M,1/c}=2A|\mathbf{a}|_{M,1/(2c)}\sum_{p=0}^{\infty}\frac{M_p}{(2c)^p}\frac{c^p}{M_p}
\\&
=2A|\mathbf{a}|_{M,1/(2c)}\sum_{p=0}^{\infty}2^{-p}=4A|\mathbf{a}|_{M,1/(2c)},
\end{align*}
where $\|\cdot\|^{\widehat{}}$ shall denote the norm on the completion $F_s$.

So we are able to define the map $T_{s}:\Lambda_{(M)}\rightarrow F_{s}$ by $\mathbf{a}\mapsto\sum_{p=0}^{\infty}a_p\chi_p$.
\qed\enddemo

Using the previous Proposition we can generalize \cite[Theorem 4.4]{surjectivity} which proves the converse implication of Theorem \ref{Beurling-Surjectivitytheorem2}.

\begin{theorem}\label{Beurling-Surjectivitytheorem4}
The inclusion
\begin{equation}\label{Beurling-Surjectivitytheorem4inclusion}
\Lambda_{(M)}\subseteq\{(f^{(rj)}(0))_{j\in\NN}: f\in\mathcal{D}_{r,(N)}([-1,1])\}
\end{equation}
implies that condition \hyperlink{SV}{$(M,N)_{SV_r}$} holds true.
\end{theorem}
Recall that inclusion \eqref{Beurling-Surjectivitytheorem4inclusion} does already imply \hyperlink{mnqr}{$(\text{nq}_r)$} for $N$.

\demo{Proof}
By Proposition \ref{Beurling-Surjectivitytheorem3} there exists a continuous linear extension mapping $T_{1}:\Lambda_{(M)}\rightarrow F_{1}$ such that $\mathcal{D}_{r,(N)}([-1,1])\ni\chi_p:=T_{1}(e^p)$ for all $p\in\NN$. By the continuity of $T_{1}$ there exists some $s\in\NN_{>0}$ and $D\ge 1$ (large enough) such that for all $\mathbf{a}=(a_p)_p\in\Lambda_{(M)}$ we get $\|T_{1}(\mathbf{a})\|^{\widehat{}}_{r,N,1}\le D|\mathbf{a}|_{M,1/s}$ (where $\|\cdot\|^{\widehat{}}$ shall again the norm in the completion as in Proposition \ref{Beurling-Surjectivitytheorem3} before). So we obtain
\begin{equation}\label{Beurling4equ0}
\forall\;p\in\NN:\;\;\;\|\chi_p\|_{r,N,1}=\|\chi_p\|^{\widehat{}}_{r,N,1}\le D|e^p|_{M,1/s}=D\frac{s^p}{M_p}.
\end{equation}
In the next step choose $1>h>0$ small enough to have $0<4hs<\frac{1}{2}$ and proceed as in the proof of Theorem \ref{Roumieu-Surjectivitytheorem4}. For $p\in\NN_{>0}$ we consider the (increasing) sequence $\tau^p:=(\tau^p_j)_{j\ge 1}$ defined by
$$\tau^p=\big(\underbrace{\left(\frac{\nu_{2p}}{h}\right)^{1/r},\dots,\left(\frac{\nu_{2p}}{h}\right)^{1/r}}_{rp-\text{times}},\underbrace{\left(\frac{\nu_{2p+1}}{h}\right)^{1/r},\dots,\left(\frac{\nu_{2p+1}}{h}\right)^{1/r}}_{r-\text{times}},\underbrace{\left(\frac{\nu_{2p+2}}{h}\right)^{1/r},\dots,\left(\frac{\nu_{2p+2}}{h}\right)^{1/r}}_{r-\text{times}},\dots\big).$$
Moreover for all $p\in\NN_{>0}$ and $0\le j\le p-1$ we define again
\begin{equation*}\label{varrhobeurling}
\varrho_{p,j}(t):=
\begin{cases}
0 &\text{if}\;t\le 0,
\\
\chi_p^{(rj)}(t)-\frac{t^{r(p-j)}}{(r(p-j))!} &\text{if}\;t>0,
\end{cases}
\end{equation*}
and consider $z>0$ satisfying
\begin{equation}\label{zbeur}
0<z<\frac{rph}{(\nu_{2p})^{1/r}}+r\sum_{l=2p+1}^{\infty}\left(\frac{h}{\nu_l}\right)^{1/r}.
\end{equation}
We use again Lemma \ref{Roumieu-Surjectivitytheorem3} for $\varrho_{p,j}$ and $a_j:=(\tau^p_j)^{-1}$ and so for each $z$ as in \eqref{zbeur} we obtain:
\begin{equation}\label{Beurling4equ1}
|\varrho_{p,j}(z)|\le\sum_{k=p}^{\infty}\frac{(4h)^k}{(\nu_{2p})^p\prod_{j=p+1}^k\nu_{p+j}}\sup_{t\in[0,z]}|\varrho^{(rk)}_{p,j}(t)|.
\end{equation}
Now estimate as follows for $k\ge p\Leftrightarrow rk\ge rp$:
\begin{align*}
\sup_{t\in[0,z]}|\varrho^{(rk)}_{p,j}(t)|&\le\sup_{t\in[0,z]}|\chi^{(r(j+k))}_{p}(t)|+1\le\|\chi_p\|_{r,N,1}N_{j+k}+1
\\&
\underbrace{\le}_{\eqref{Beurling4equ0}}D|e^p|_{M,1/s}N_{j+k}+1=D s^p\frac{N_{j+k}}{M_p}+1\le D s^k\frac{N_{j+k}}{M_p}+1\le 2D s^k\frac{N_{j+k}}{M_p},
\end{align*}
where in the last estimate we have used $M_p\le N_p\le N_{j+k}$. As in the Roumieu case we are going to prove for each $z$ as in \eqref{zbeur}
\begin{equation*}\label{Beurling4equ3}
|\varrho_{p,j}(z)|\le 2D\frac{N_j}{M_p}2^{-p+1},
\end{equation*}
which holds by
\begin{align*}
|\varrho_{p,j}(z)|&\underbrace{\le}_{\eqref{Beurling4equ1},\eqref{Roumieutheorem4equ2}}\sum_{k=p}^{\infty}(4h)^k\frac{N_j}{N_{j+k}}\sup_{t\in[0,z]}|\varrho^{(rk)}_{p,j}(t)|\le\sum_{k=p}^{\infty}(4h)^k\frac{N_j}{N_{j+k}}2Ds^{k}\frac{N_{j+k}}{M_p}
\\&
\le 2D\frac{N_j}{M_p}\sum_{k=p}^{\infty}(\underbrace{4hs}_{<1/2})^k\le 2D\frac{N_j}{M_p}2^{-p+1}.
\end{align*}
For any $z>0$ as in \eqref{zbeur} we are going to prove as in the Roumieu case
\begin{equation}\label{Beurling4equ4}
\forall\;p\in\NN_{>0}\;\forall\;j\in\NN, 0\le j\le p-1:\;\;\frac{1}{2}\frac{z^{r(p-j)}}{(r(p-j))!}\le 2D\frac{s^p N_j}{M_p}.
\end{equation}
Case $1$ is completely the same as above, for case $2$ we have that $\chi_p^{(rj)}(z)>\frac{1}{2}\frac{z^{r(p-j)}}{(r(p-j))!}$ implies
\begin{align*}
\frac{1}{2}\frac{z^{r(p-j)}}{(r(p-j))!}&<\chi_p^{(rj)}(z)\le\|\chi_p\|_{r,N,1}N_j\le D|e^p|_{M,1/s}N_j=D\frac{s^p N_j}{M_p},
\end{align*}
where we have used \eqref{Beurling4equ0}.\vspace{6pt}

Using \eqref{Beurling4equ4} for the choice $z:=\sum_{l=2p+1}^{\infty}\left(\frac{h}{\nu_l}\right)^{1/r}$ yields the conclusion by the same proof as in the Roumieu case above.
\qed\enddemo

\section{Special cases and consequences}\label{section6}
\subsection{The constant case $M=N$}\label{section61}
We are going to apply the results from Sections \ref{section4} and \ref{section5} to the case $M=N$. In this case \hyperlink{gammar}{$(M,N)_{\gamma_r}$} is precisely condition \hyperlink{gammar}{$(\gamma_r)$} from \cite{Schmetsvaldivia00} (see Section \ref{section21}) and implies \hyperlink{SV}{$(M,M)_{SV_r}$}. The special case $r=1$ yields \hyperlink{gamma1}{$(\gamma_1)$} from \cite{petzsche} respectively $(\ast)$ in \cite[p. 385]{surjectivity} with $M^1=M^2=M$.\vspace{6pt}

So we can reformulate and generalize \cite[Theorem 3.6]{surjectivity} where only the Roumieu case was considered (for $r=1$) and even in this situation we have a slightly more general statement because \hyperlink{mnq}{$(\text{nq})$} on $M$ is not assumed in our result. As mentioned in \cite{surjectivity} the case $r=1$ is reproving one of the main results from \cite{petzsche}. But also there assumption \hyperlink{mnq}{$(\text{nq})$} on $M$, which was called $(\gamma)$ and a basic property, is superfluous as we have already commented in \cite[Def. 4.3, Thm. 4.4, p. 154]{injsurj}.

Note that the following result for $r\ge 2$ provides a characterization of the surjectivity of the restriction mapping in terms of condition \hyperlink{gammar}{$(\gamma_r)$} which completes the results obtained for the Roumieu case in \cite{Schmetsvaldivia00}.\vspace{6pt}


Let $M\in\RR_{>0}^{\NN}$ and $r\in\NN_{\ge 1}$ satisfying

\begin{itemize}
\item[$(I)$] $M\in\hyperlink{LCset}{\mathcal{LC}}$,


\item[$(II)_{B,r}$] $\lim_{p\rightarrow\infty}\left(\frac{(M_p)^{1/r}}{p!}\right)^{1/p}=\infty$.
\end{itemize}

\begin{theorem}\label{specialcase1}
Let $M$ be as assumed above and $r\in\NN_{\ge 1}$, then the following are equivalent:

\begin{itemize}
\item[$(i)$] The Borel map $\mathcal{B}^r:\mathcal{D}_{r,\{M\}}([-1,1])\rightarrow\Lambda_{\{M\}}$, $f\mapsto(f^{(rj)}(0))_{j\in\NN}$, is surjective.

\item[$(ii)$] The Borel map $\mathcal{B}^r:\mathcal{D}_{r,(M)}([-1,1])\rightarrow\Lambda_{(M)}$, $f\mapsto(f^{(rj)}(0))_{j\in\NN}$, is surjective.

\item[$(iii)$] Condition \hyperlink{SV}{$(M,M)_{SV_r}$} holds true.

\item[$(iv)$] Condition \hyperlink{gammarmix}{$(M,M)_{\gamma_r}$}, i.e. \hyperlink{gammar}{$(\gamma_r)$} from \cite{Schmetsvaldivia00}, holds true for $M$.
\end{itemize}
\end{theorem}

For treating the Roumieu case it is sufficient to replace $(II)_{B,r}$ by the weaker assumption $(II)_{R,r}$, i.e. $\liminf_{p\rightarrow\infty}\left(\frac{(M_p)^{1/r}}{p!}\right)^{1/p}>0$, and then in Theorem \ref{specialcase1} we obtain $(i)\Leftrightarrow(iii)\Leftrightarrow(iv)$.

\demo{Proof}
By Theorems \ref{Roumieu-Surjectivitytheorem}, \ref{Beurling-Surjectivitytheorem} and Lemma \ref{lambdaexpression} it remains to prove $(iii)\Rightarrow(iv)$.

Put $\lambda_{p,s}:=\lambda^{M,M}_{p,s}$ and for all $s,p\in\NN_{>0}$ we get
\begin{align*}
\lambda_{2p,s}&\underbrace{\ge}_{j=p}\left(\frac{M_{2p}}{s^{2p} M_p}\right)^{1/p}=s^{-2}(\mu_{p+1}\cdots\mu_{2p})^{1/p}\ge s^{-2}(\mu_p)^{p/p}=s^{-2}\mu_p,
\end{align*}
hence we estimate for any $p\in\NN_{>0}$ as follows:
\begin{align*}
\frac{(\mu_p)^{1/r}}{p}\sum_{k=p}^{\infty}\left(\frac{1}{\mu_k}\right)^{1/r}&=\frac{(\mu_p)^{1/r}}{p}\underbrace{\sum_{k=p}^{2p-1}\left(\frac{1}{\mu_k}\right)^{1/r}}_{\le p/(\mu_p)^{1/r}}+\frac{(\mu_p)^{1/r}}{p}\sum_{k=2p}^{\infty}\left(\frac{1}{\mu_k}\right)^{1/r}
\\&
\le 1+s^{2/r}\frac{(\lambda_{2p,s})^{1/r}}{p}\sum_{k=2p}^{\infty}\left(\frac{1}{\mu_k}\right)^{1/r}\le 1+s^{2/r}2D,
\end{align*}
where the last inequality holds for some $D\ge 1$ large. This proves \hyperlink{gammar}{$(\gamma_r)$} for $M$.
\qed\enddemo

\subsection{More ultradifferentiable $r$-ramification function spaces}\label{testfunctionspaces}
So far we have considered the mixed Borel mapping setting using the classes $\mathcal{D}_{r,\{N\}}([-1,1])$ and $\mathcal{D}_{r,(N)}([-1,1])$. But in \cite[Section 3]{Schmetsvaldivia00} also the following closely related spaces have been introduced in order to prove extension theorems in the ultraholomorphic setting. Moreover, in \cite[Section 4]{injsurj} the authors have used these methods and spaces to treat the surjectivity of the asymptotic Borel map by giving a connection between \hyperlink{gammar}{$(\gamma_r)$} and the growth index $\gamma(M)$ introduced in \cite[Def. 1.3.5]{Thilliezdivision}.

We put
\begin{align*}
\mathcal{L}_{r,\{M\}}:=&\{f\in\mathcal{E}([0,\infty),\CC):\;\supp(f)\subseteq[0,1],\; f^{(rn+j)}(0)=0,\;\forall\;n\in\NN\;\forall\;j\in\{1,\dots,r-1\},
\\&
\exists\;h>0:\;\sup_{n\in\NN, x\in[0,\infty)}\frac{|f^{(rn)}(x)|}{h^nM_n}<\infty\},
\end{align*}
\begin{align*}
\mathcal{N}_{r,\{M\}}:=&\{f\in\mathcal{E}([0,\infty),\CC):\;f^{(rn+j)}(0)=0,\;\forall\;n\in\NN\;\forall\;j\in\{1,\dots,r-1\},
\\&
\exists\;h>0:\;\sup_{n\in\NN, x\in[0,\infty)}\frac{|f^{(rn)}(x)|}{h^nM_n}<\infty\},
\end{align*}
\begin{align*}
\mathcal{E}_{r,\{M\}}:=&\{f\in\mathcal{E}([0,1],\CC):\;f^{(rn+j)}(0)=0,\;\forall\;n\in\NN\;\forall\;j\in\{1,\dots,r-1\},
\\&
\exists\;h>0:\;\sup_{n\in\NN, x\in[0,1]}\frac{|f^{(rn)}(x)|}{h^nM_n}<\infty\},
\end{align*}
respectively for the Beurling type classes where $\exists\;h>0$ is replaced by $\forall\;h>0$.

We will now see that Theorem \ref{Roumieu-Surjectivitytheorem}, respectively Theorem \ref{Beurling-Surjectivitytheorem}, remains true if we replace $\mathcal{D}_{r,\{N\}}([-1,1])$ by any of the new classes above (respectively for the Beurling case).\vspace{6pt}

First, we note that
\begin{equation}\label{spaceinclusion}
\mathcal{D}_{r,[M]}([-1,1])\subseteq\mathcal{L}_{r,[M]}\subseteq\mathcal{N}_{r,[M]}\subseteq\mathcal{E}_{r,[M]}\subseteq\mathcal{E}_{[P^{M,r}]}([0,1],\CC).
\end{equation}
The first three inclusions follow immediately by definition and restriction. The last inclusion was shown in \cite[Prop. 5.2]{Schmetsvaldivia00} for the Roumieu, and in \cite[Prop. 4.2]{Schmetsvaldivia00} for the Beurling case by using the so-called {\itshape Gorny-Cartan-inequalities}, e.g. see \cite[6.4. IV]{mandelbrojtbook}.\vspace{6pt}

An immediate consequence of the first inclusion is that Theorem \ref{Roumieu-Surjectivitytheorem1} in the Roumieu, respectively Theorem \ref{Beurling-Surjectivitytheorem2} in the Beurling case, can be generalized as follows:

\begin{theorem}\label{testfunctionspacesthm1}
Let $M$ and $N$ be given satisfying $(I),(III)$ and $(II)_{R,r}$ in the Roumieu, or $(I),(III)$ and $(II)_{B,r}$ in the Beurling case for some $r\in\NN_{\ge 1}$. If \hyperlink{SV}{$(M,N)_{SV_r}$} holds true, then we get
$$\Lambda_{[M]}\subseteq\{(f^{(rj)}(0))_{j\in\NN}: f\in\mathcal{A}\},$$
with $\mathcal{A}\in\{\mathcal{L}_{r,[N]},\mathcal{N}_{r,[N]},\mathcal{E}_{r,[N]}\}$.
\end{theorem}

On the other hand, all the Roumieu type classes above are countable $(LB)$-spaces and the Beurling type classes are Fréchet spaces, hence Proposition \ref{Roumieu-Surjectivitytheorem2} in the Roumieu and Proposition \ref{Beurling-Surjectivitytheorem3} in the Beurling case remain true if we replace $\mathcal{D}_{r,[N]}([-1,1])$ by any of the spaces above (of the particular type). And also Theorem \ref{Roumieu-Surjectivitytheorem4} in the Roumieu case and Theorem \ref{Beurling-Surjectivitytheorem4} in the Beurling case will follow with the same proofs.\vspace{6pt}

An immediate consequence of these results is that the inclusion $\Lambda_{[M]}\subseteq\{(f^{(rj)}(0))_{j\in\NN}: f\in\mathcal{A}\}$, $\mathcal{A}\in\{\mathcal{L}_{r,[N]},\mathcal{N}_{r,[N]},\mathcal{E}_{r,[N]}\}$, implies \hyperlink{SV}{$(M,N)_{SV_r}$} and as already seen in $(iii)$ in Remark \ref{Roumieu-Surjectivitytheoremremark} we get condition \hyperlink{mnqr}{$(\on{nq}_r)$} for $N$. So in any case \hyperlink{mnqr}{$(\on{nq}_r)$} is not required as a standard assumption for the weight $N$.\vspace{6pt}

Thus we can summarize all our results in the following final statement:

\begin{theorem}\label{finaltheorem}
Let $M,N\in\RR_{>0}^{\NN}$ be given satisfying $(I),(III)$ and $(II)_{R,r}$ in the Roumieu, or $(I),(III)$ and $(II)_{B,r}$ in the Beurling case for some $r\in\NN_{\ge 1}$. Then the following are equivalent:

\begin{itemize}
\item[$(i)$] $\Lambda_{[M]}\subseteq\{(f^{(rj)}(0))_{j\in\NN}: f\in\mathcal{A}\}$, $\mathcal{A}\in\{\mathcal{D}_{r,[N]}([-1,1]),\mathcal{L}_{r,[N]},\mathcal{N}_{r,[N]},\mathcal{E}_{r,[N]}\}$,

\item[$(ii)$] condition \hyperlink{SV}{$(M,N)_{SV_r}$} does hold true.
\end{itemize}
\end{theorem}

We point out that the special case $M=N$ yields the characterization of the surjectivity of the Borel map for all $r$-ramification spaces in terms of condition \hyperlink{gammar}{$(\gamma_r)$}, what for some of these spaces completes the results in \cite{Schmetsvaldivia00} and/or weakens their assumptions.

\subsection{An application to a Whitney extension theorem in the mixed setting}\label{section62}
Using the ramification conditions in this present work we are able now to reformulate and generalize the results from \cite[Section 5.4]{whitneyextensionmixedweightfunction} and in this section the restriction $r\in\NN_{\ge 1}$ will be not necessary since no ramification spaces are involved. First we have to introduce the notion of an associated weight function.

Let $M\in\RR_{>0}^{\NN}$ (with $M_0=1$), then the {\itshape associated function} $\omega_M: \RR_{\ge 0}\rightarrow\RR\cup\{+\infty\}$ is defined by
\begin{equation*}\label{assofunc}
\omega_M(t):=\sup_{p\in\NN}\log\left(\frac{t^p}{M_p}\right)\;\;\;\text{for}\;t>0,\hspace{30pt}\omega_M(0):=0.
\end{equation*}
For an abstract introduction of the associated function we refer to \cite[Chapitre I]{mandelbrojtbook}, see also \cite[Definition 3.1]{Komatsu73}. If $\liminf_{p\rightarrow\infty}(M_p)^{1/p}>0$, then $\omega_M(t)=0$ for sufficiently small $t$, since $\log\left(\frac{t^p}{M_p}\right)<0\Leftrightarrow t<(M_p)^{1/p}$ holds for all $p\in\NN_{>0}$. Moreover under this assumption $t\mapsto\omega_M(t)$ is a continuous nondecreasing function, which is convex in the variable $\log(t)$ and tends faster to infinity than any $\log(t^p)$, $p\ge 1$, as $t\rightarrow\infty$. $\lim_{p\rightarrow\infty}(M_p)^{1/p}=\infty$ implies that $\omega_M(t)<\infty$ for each finite $t$, and this shall be considered as a basic assumption for defining $\omega_M$. One may also introduce the function $\Sigma_M:\RR_{\ge 0}\rightarrow\RR$ defined by
$$\Sigma_M(t):=|\{p\in\NN_{>0}: \mu_p\le t\}|,\quad \textrm{($|\ |$ denotes cardinal)}$$
which allows us to write
\begin{equation}\label{eq:OmNuineq}
  \omega_M(t)=\int_{0}^{t}\Sigma_M(u)\frac{du}{u},\quad t>0.
\end{equation}

For all $t,s>0$ we get
\begin{equation}\label{omegaMspower}
\omega_M(t^s)=\sup_{p\in\NN}\log\left(\frac{t^{sp}}{M_p}\right)=\sup_{p\in\NN}\log\left(\left(\frac{t^p}{(M_p)^{1/s}}\right)^s\right)=s\omega_{M^{1/s}}(t),
\end{equation}
with $M^{1/s}_p=(M_p)^{1/s}$.

We can generalize \cite[Lemma 5.7]{whitneyextensionmixedweightfunction} as follows.

\begin{lemma}\label{mglemma1}
Let $M,N\in\hyperlink{LCset}{\mathcal{LC}}$ be given with $\mu\le\nu$ (which implies $M\le N$) and satisfying \hyperlink{gammarmix}{$(M,N)_{\gamma_r}$} for some arbitrary $r>0$. Then the associated weight functions are satisfying
\begin{equation}\label{mglemma1equ}
\exists\;C>0\;\forall\;t\ge 0:\;\;\;\int_1^{\infty}\frac{\omega_N(tu)}{u^{1+1/r}}du\le C\omega_M(t)+C.
\end{equation}
\end{lemma}

\demo{Proof}
By definition \hyperlink{gammarmix}{$(M,N)_{\gamma_r}$} holds true if and only if \hyperlink{gammarmix}{$(M^{1/r},N^{1/r})_{\gamma_1}$} holds, which is precisely \cite[(5.17)]{whitneyextensionmixedweightfunction}. Since $M^{1/r},N^{1/r}\in\hyperlink{LCset}{\mathcal{LC}}$ and $\mu^{1/r}\le\nu^{1/r}$ (with $\mu_p^{1/r}=(\mu_p)^{1/r}$), by \cite[Lemma 5.7]{whitneyextensionmixedweightfunction} we get
\begin{equation*}\label{mglemma1equ1}
\exists\;C>0\;\forall\;t\ge 0:\;\;\;\int_1^{\infty}\frac{\omega_{N^{1/r}}(tu)}{u^2}du\le C\omega_{M^{1/r}}(t)+C.
\end{equation*}
By applying \eqref{omegaMspower} the right-hand side gives $\frac{C}{r}\omega_M(t^r)+C$, whereas the left-hand side gives $$\int_1^{\infty}\frac{\omega_{N^{1/r}}(tu)}{u^2}du=\frac{1}{r}\int_1^{\infty}\frac{\omega_{N}((tu)^r)}{u^2}du=\frac{1}{r^2}\int_1^{\infty}\frac{\omega_{N}(t^rv)}{v^{1+1/r}}dv.$$ Finally, by replacing on both sides $t^r$ by $t$, we are done.
\qed\enddemo

Now we are proving the converse statement, here we have to make use of \hyperlink{mg}{$(\on{mg})$} and we are generalizing \cite[Corollary 4.6 $(iii)$]{index} to a mixed setting.

\begin{lemma}\label{mglemma2}
Let $M,N\in\hyperlink{LCset}{\mathcal{LC}}$ be given with $\mu\le\nu$ (which implies $M\le N$ and $\mu^r\le\nu^r$ for all $r>0$), such that $M$ does have \hyperlink{mg}{$(\on{mg})$} and, moreover, the condition~\eqref{mglemma1equ} is satisfied. Then, \hyperlink{gammarmix}{$(M,N)_{\gamma_r}$} holds.
\end{lemma}

\demo{Proof}
To avoid technical complications in the proof we assume that the sequences $p\mapsto\mu_p$ and $p\mapsto\nu_p$ are strictly increasing. This can be done w.l.o.g. by passing, if necessary, to equivalent sequences as follows: First one can construct $\widetilde{\mu}=(\widetilde{\mu}_p)_p$ and $\widetilde{\nu}=(\widetilde{\nu}_p)_p$ such that $\mu_p\le\widetilde{\mu}_p\le a\mu_p$ and $\nu_p\le\widetilde{\nu}_p\le a\nu_p$ holds true for some constant $a>1$ and all $p\in\NN$ and both $p\mapsto\widetilde{\mu}_p$ and $p\mapsto\widetilde{\nu}_p$ are strictly increasing. Then put $\widehat{\mu}_p:=a^{-1}\widetilde{\mu}_p$, $\widehat{\nu}_p:=a\widetilde{\nu}_p$, hence both $p\mapsto\widehat{\mu}_p$, $p\mapsto\widehat{\nu}_p$ are still strictly increasing, and $\mu\le\nu$ implies $\widehat{\mu}\le\widetilde{\mu}\le a\widetilde{\nu}\le\widehat{\nu}$. By considering $\widehat{M}$ and $\widehat{N}$ in the obvious way, we see that $\widehat{M}$ is equivalent to $M$, $\widehat{N}$ is equivalent to $N$ and also the inequality~\eqref{mglemma1equ} is preserved with $\omega_{M}$ substituted by $\omega_{\widehat{M}}$ and $\omega_{N}$ substituted by $\omega_{\widehat{N}}$, since $\widehat{M}\le M$, $\widehat{N}\ge N$, and so $\omega_{\widehat{M}}(t)\ge\omega_{M}(t)$ and $\omega_{\widehat{N}}(t)\le \omega_{N}(t)$ for every $t>0$. If we can deduce from this inequality that \hyperlink{gammarmix}{$(\widehat{M},\widehat{N})_{\gamma_r}$} holds, it is clear that also \hyperlink{gammarmix}{$(M,N)_{\gamma_r}$} will be valid, as desired.\vspace{6pt}

Since $M$ does have \hyperlink{mg}{$(\on{mg})$} (which is preserved by switching to any equivalent sequence), by the estimate given in the proof of \cite[Theorem 4.4 $(ii)$]{index}, we have $\omega_M(t)\le A\Sigma_M(t)+A$ for some $A\ge 1$ and all $t\ge 0$. Hence, replacing $t$ by $t^r$ in~\eqref{mglemma1equ}, we get
\begin{equation}\label{eq_mixed_omegaMN_mg}
\int_1^{\infty}\frac{\omega_N(t^ru)}{u^{1+1/r}}du\le C\omega_M(t^r)+C\le C_1\Sigma_M(t^r)+C_1
\end{equation}
for some $C_1>0$ (large) and all $t\ge 0$.
The monotonicity of $\Sigma_N$ and \eqref{eq:OmNuineq} together imply
\begin{equation*}
  \omega_N(et)=\int_{0}^{et}\Sigma_N(u)\frac{du}{u}\geq \int_{t}^{et}\Sigma_N(u)\frac{du}{u}\ge \Sigma_N(t).
\end{equation*}
So, the left-hand side in~\eqref{eq_mixed_omegaMN_mg} can be estimated as
\begin{equation}\label{eq_omegaN_sigmaN}
\int_1^{\infty}\frac{\omega_N(t^ru)}{u^{1+1/r}}du\ge
\int_1^{\infty}\frac{\Sigma_N(t^ru/e)}{u^{1+1/r}}du= \frac{t}{e^{1/r}}\int_{t^r/e}^{\infty}\frac{\Sigma_N(v)}{v^{1+1/r}}dv\ge
\frac{t}{e^{1/r}}\int_{t^r}^{\infty}\frac{\Sigma_N(v)}{v^{1+1/r}}dv.
\end{equation}



If $t\ge \nu_1^{1/r}$ and we put $p=\Sigma_N(t^r)\ge 1$, we may compute and estimate the last integral as
\begin{align*}
&\int_{t^r}^{\infty}\frac{\Sigma_N(v)}{v^{1+1/r}}dv= \int_{t^r}^{\nu_{p+1}}\frac{\Sigma_N(v)}{v^{1+1/r}}dv+ \sum_{j=p}^{\infty}\int_{\nu_{j+1}}^{\nu_{j+2}} \frac{\Sigma_N(v)}{v^{1+1/r}}dv
\\&
=p\left[-r\frac{1}{v^{1/r}}\right]_{v=t^r}^{v=\nu_{p+1}}+ \sum_{j=p}^{\infty}(j+1) \left[-r\frac{1}{v^{1/r}}\right]_{v=\nu_{j+1}}^{v=\nu_{j+2}}
= \frac{pr}{t}+ r\sum_{j=p}^{\infty}\left(\frac{1}{\nu_{j+1}}\right)^{1/r}\\
&\ge r\sum_{j=\Sigma_N(t^r)+1}^{\infty}\left(\frac{1}{\nu_{j}}\right)^{1/r}.
\end{align*}
Gathering this with \eqref{eq_mixed_omegaMN_mg} and \eqref{eq_omegaN_sigmaN} we deduce that
$$
\frac{rt}{e^{1/r}}\sum_{j=\Sigma_N(t^r)+1}^{\infty} \left(\frac{1}{\nu_{j}}\right)^{1/r}\le C_1\Sigma_M(t^r)+C_1.
$$
Since $\mu\le\nu$ we have $\Sigma_M(t)\ge\Sigma_N(t)$ for every $t$, and so for all $t\ge \nu_1^{1/r}$ we have the inequality
\begin{equation}\label{eq_ineq_gammaMN_SigmaMN}
\frac{t}{\Sigma_M(t^r)}\sum_{j=\Sigma_M(t^r)+1}^{\infty} \left(\frac{1}{\nu_{j}}\right)^{1/r}\le
\frac{t}{\Sigma_M(t^r)}\sum_{j=\Sigma_N(t^r)+1}^{\infty} \left(\frac{1}{\nu_{j}}\right)^{1/r}\le
2C_1\frac{e^{1/r}}{r}.
\end{equation}
For any $q\in\NN_{>0}$, we choose $t=\mu_q^{1/r}$ in~\eqref{eq_ineq_gammaMN_SigmaMN} and deduce that
$$
\frac{\mu_q^{1/r}}{q}\sum_{j=q}^\infty\frac{1}{\nu_j^{1/r}}\le
1+\frac{\mu_q^{1/r}}{q}\sum_{j=q+1}^\infty\frac{1}{\nu_j^{1/r}}
\le 1+2C_1\frac{e^{1/r}}{r},
$$
as desired.

\qed\enddemo

The next result characterizes the possibility of obtaining mixed Whitney extension results in terms of the mixed conditions with a ramification parameter $r$. It generalizes to any $r>0$ the result for $r=1$ obtained by A. Rainer and the third author \cite[Theorem 5.9]{whitneyextensionmixedweightfunction}, and it improves it by dropping the moderate growth condition for $N$. Also, for $r=1$ one recovers the central theorem in the Roumieu version of \cite{ChaumatChollet94}, which is indeed used in our arguments.
We are considering, for a compact set $E\subseteq\RR^n$, the class $\mathcal{B}_{\{M\}}(E)$ of {\itshape Whitney ultrajets} of Roumieu type defined by $M$, for a precise definition we refer to \cite[Definition 2.7]{whitneyextensionmixedweightfunction}. Finally let $j^{\infty}_E$ be the {\itshape jet mapping} which assigns to each smooth function $f$ defined in $\RR^n$ the infinite jet consisting of its partial derivatives of all orders restricted to $E$ (i.e. $j^{\infty}_{\{x_0\}}$ is the Borel map $\mathcal{B}_{x_0}$).

\begin{theorem}\label{mgtheorem1}
Let $M,N\in\hyperlink{LCset}{\mathcal{LC}}$ be given with $\mu\le\nu$ and such that $M$ satisfies \hyperlink{mg}{$(\on{mg}$)}. Moreover assume that
\begin{equation}\label{mgtheorem1equ1}
\exists\;r>0\;\exists\;C\ge 1\;\forall\;1\le j\le k:\:\:\:\frac{(\mu_j)^{1/r}}{j}\le C\frac{(\mu_k)^{1/r}}{k},\hspace{20pt}\frac{(\nu_j)^{1/r}}{j}\le C\frac{(\nu_k)^{1/r}}{k}
\end{equation}
and that $\lim_{j\rightarrow\infty}\frac{(\mu_j)^{1/r}}{j}=\infty$.

Then the following conditions are equivalent ($r>0$ denoting the number arising in \eqref{mgtheorem1equ1}):

\begin{itemize}
\item[$(i)$] For every compact set $E\subseteq\RR^n$ we get $j^{\infty}_E(\mathcal{B}_{\{N^{1/r}\}}(\RR^n))\supseteq\mathcal{B}_{\{M^{1/r}\}}(E)$.

\item[$(ii)$] The associated weight functions satisfy $\exists\;C>0\;\forall\;t\ge 0:\;\;\;\int_1^{\infty}\frac{\omega_N(tu)}{u^{1+1/r}}du\le C\omega_M(t)+C$.

\item[$(iii)$] Condition \hyperlink{gammarmix}{$(M,N)_{\gamma_r}$} holds true.

\item[$(iv)$] Condition \hyperlink{SV}{$(M,N)_{SV_r}$} holds true.
\end{itemize}
If in the assumption above $r\in\NN_{\ge 1}$ and $P^{M,r}$, $P^{N,r}$ are denoting the corresponding $r$-interpolating sequences, then moreover $(i)-(iv)$ are equivalent to
\begin{itemize}
\item[$(i')$] For every compact set $E\subseteq\RR^n$ we get $j^{\infty}_E(\mathcal{B}_{\{P^{N,r}\}}(\RR^n))\supseteq\mathcal{B}_{\{P^{M,r}\}}(E)$.

\item[$(ii')$] The associated weight functions satisfy $\exists\;C>0\;\forall\;t\ge 0:\;\;\;\int_1^{\infty}\frac{\omega_{P^{N,r}}(tu)}{u^{2}}du\le C\omega_{P^{M,r}}(t)+C$.

\item[$(iii')$] Condition \hyperlink{gammarmix}{$(P^{M,r},P^{N,r})_{\gamma_1}$} holds true.

\item[$(iv')$] Condition \hyperlink{SV}{$(P^{M,r},P^{N,r})_{SV_1}$} holds true.
\end{itemize}
\end{theorem}

\begin{remark}\label{mgtheorem1remark}
\eqref{mgtheorem1equ1} does precisely mean that both sequences $(\mu_j/j^r)_{j\ge 1}$ and $(\nu_j/j^r)_{j\ge 1}$ are almost increasing. As shown in \cite[Theorem 3.11]{index} we have $\gamma(M),\gamma(N)\ge r$ with $\gamma(M)$ denoting the growth index introduced in \cite{Thilliezdivision}, see also \cite{injsurj}, and moreover we can replace $M$ and $N$ by equivalent sequences $\widetilde{M}$ and $\widetilde{N}$ such that $j\mapsto\frac{\widetilde{\mu_j}}{j^r}$, $j\mapsto\frac{\widetilde{\nu_j}}{j^r}$ are nondecreasing, i.e. $((\widetilde{M}_j)^{1/r}/j!)_j$ and $((\widetilde{N}_j)^{1/r}/j!)_j$ are log-convex.
\end{remark}

\demo{Proof}
First, as seen in Section \ref{section32}, whenever $M$ has \hyperlink{mg}{$(\on{mg})$} we get the equivalence of $(iii)$ and $(iv)$.
Moreover, Lemmas~\ref{mglemma1} and~\ref{mglemma2} show the equivalence of $(ii)$ and $(iii)$ under the same condition.

Assume now that \hyperlink{gammarmix}{$(M,N)_{\gamma_r}$} holds true. Then \hyperlink{mnqr}{$(\text{nq}_r)$} holds for $N$ and by Remark \ref{mgtheorem1remark} and assumption \eqref{mgtheorem1equ1} we can replace $M$ and $N$ by $\widetilde{M}$ and $\widetilde{N}$ such that $((\widetilde{M}_j)^{1/r}/j!)_j$ and $((\widetilde{N}_j)^{1/r}/j!)_j$ are log-convex. Equivalence preserves \hyperlink{mg}{$(\on{mg}$)} for $\widetilde{M}$, \hyperlink{mnqr}{$(\text{nq}_r)$} for $\widetilde{N}$ (by Carleman's inequality) and finally \hyperlink{gammarmix}{$(\widetilde{M},\widetilde{N})_{\gamma_r}$}. Thus we are able to apply \cite[Theorem 11]{ChaumatChollet94} to $M\equiv((\widetilde{M}_j)^{1/r}/j!)_j$ and $M'\equiv((\widetilde{N}_j)^{1/r}/j!)_j$ (our notation for weight sequences differs from the one used in \cite{ChaumatChollet94} by a factorial term), which yields $(i)$.
\vspace{6pt}

For proving $(i)\Rightarrow(iii)$ we just need to apply \cite[Proposition 27]{ChaumatChollet94} again to $M\equiv((\widetilde{M}_j)^{1/r}/j!)_j$ and $M'\equiv((\widetilde{N}_j)^{1/r}/j!)_j$ with $E=\{0\}$.

Please note that assumption \hyperlink{mnqr}{$(\text{nq}_r)$} on $N$ (i.e. \hyperlink{mnq}{$(\text{nq})$} for $N^{1/r}$ and $\widetilde{N}^{1/r}$, as it was required in \cite[Proposition 27]{ChaumatChollet94}) is not necessary in our statement (a similar argument can be found explicitly in \cite[Corollary 2]{whitneyextensionmixedweightfunctionII} and implicitly in \cite[Theorem 5.9 $(1)\Rightarrow(3)$]{whitneyextensionmixedweightfunction}): First, if $M^{1/r}$ does satisfy \hyperlink{mnq}{$(\text{nq})$}, then by $\mu\le\nu\Leftrightarrow\mu^{1/r}\le\nu^{1/r}$ also $N^{1/r}$ shares this property (or equivalently $N$ has \hyperlink{mnqr}{$(\text{nq}_r)$}).

Second, suppose $M^{1/r}$ does not have \hyperlink{mnq}{$(\text{nq})$} and assume that also $N^{1/r}$ does not have \hyperlink{mnq}{$(\text{nq})$}. In this case $\lim_{j\rightarrow\infty}\frac{(\mu_j)^{1/r}}{j}=\infty$, which implies $\lim_{j\rightarrow\infty}(M^{1/r}_j/j!)^{1/j}=\infty$ by \eqref{mgequiv} and Stirling's formula, together with the main result \cite[Theorem 2]{borelmappingquasianalytic} applied to $M^{1/r}$ and $N^{1/r}$ (instead of $M$ and $N$) yield a contradiction to assumption $(i)$ applied to the case $E=\{0\}$.

Note that by $\lim_{j\rightarrow\infty}\frac{(\mu_j)^{1/r}}{j}=\infty$ the (local) ultradifferentiable class defined by $M^{1/r}$ cannot coincide with the real-analytic functions.
\vspace{6pt}

For the additional part we remark that by Lemma \ref{prmg} also both $P^{M,r},P^{N,r}\in\hyperlink{LCset}{\mathcal{LC}}$ do have \hyperlink{mg}{$(\on{mg}$)}, $\pi^{M,r}\le\pi^{N,r}$ holds by $\mu\le\nu$ and \eqref{interpolatingsequ0}. Since $(P^{M,r}_j)^{1/j}\le\pi^{M,r}_j\le A(P^{M,r}_j)^{1/j}$ for some $A\ge 1$ and all $j\in\NN_{\ge 1}$, by $(iii)$ in Remark \ref{beurlingcondtionremark}, $\lim_{j\rightarrow\infty}\frac{(\mu_j)^{1/r}}{j}=\infty$ does imply $\lim_{j\rightarrow\infty}\frac{\pi^{M,r}_j}{j}=\infty$ (and for $P^{N,r}$ too). Finally, concerning \eqref{mgtheorem1}, for given $1\le l_1\le l_2$ we have $l_1=rk_1+j$, $l_2=rk_2+j$ for some $k_1,k_2\in\NN$ with $k_1\le k_2$ and $j\in\{1,\dots,r\}$ and use \eqref{interpolatingsequ0} to estimate by
$$\frac{\pi^{M,r}_{l_1}}{l_1}=\frac{(\mu_{k_1+1})^{1/r}}{rk_1+j}\le\frac{(\mu_{k_1+1})^{1/r}}{k_1+1}\le C\frac{(\mu_{k_2+1})^{1/r}}{k_2+1}\le rC\frac{(\mu_{k_2+1})^{1/r}}{rk_2+j}=rC\frac{\pi^{M,r}_{l_2}}{l_2},$$
because $\frac{rk_2+j}{k_2+1}\le\frac{rk_2+r}{k_2+1}=r$ for all $k_2\in\NN$ and $j\in\{1,\dots,r\}$.

Hence we can apply the above arguments to $P^{M,r}$ and $P^{N,r}$ instead of $M$ and $N$ and the equivalences $(i')-(iv')$ follow. Finally, as shown in Lemma \ref{prstrongnonquasi}, we have $(iii)\Leftrightarrow(iii')$.
\qed\enddemo


%



\section{Characterization of (non)quasianalyticity for $r$-ramification ultradifferentiable classes}\label{nonquasi}
The aim of this final section is to characterize the nonquasianalyticity resp. nontriviality of all classes of ultradifferentiable functions defined in \eqref{Roumieu-LB-space} and \eqref{Beurling-Frechet-space} and Section \ref{testfunctionspaces}. As mentioned in Remark \ref{remarktestfunctionnontrivial} we could equivalently replace in the main result below $\mathcal{D}_{r,[M]}([-1,1])$ by $\mathcal{D}_{r,[M]}([-a,a])$ for arbitrary $a>0$ and simultaneously replace $1$ by $a$ in the definition of $\mathcal{E}_{r,[M]}$ and $\mathcal{L}_{r,[M]}$.

\begin{theorem}\label{nonquasitheorem}
Let $M\in\hyperlink{LCset}{\mathcal{LC}}$ and $r\in\NN_{\ge 1}$. Then the classes $\mathcal{D}_{r,[M]},\mathcal{L}_{r,[M]},\mathcal{N}_{r,[M]}$ and $\mathcal{E}_{r,[M]}$ are nonquasianalytic, i.e. the Borel map restricted to any of these classes is not injective, if and only if $M$ satisfies \hyperlink{mnqr}{$(\on{nq}_r)$}.
\end{theorem}

\demo{Proof}
Suppose $M$ satisfies \hyperlink{mnqr}{$(\on{nq}_r)$}. By the inclusions from \eqref{spaceinclusion}, it suffices to prove the nonquasianalyticity for $\mathcal{D}_{r,[M]}([-1,1])$. And this fact has already been shown in Lemma \ref{Petzsche2remark} (recall that in this situation $P^{M,r}$ satisfies \hyperlink{mnq}{$(\on{nq})$} and so $\mathcal{E}_{[P^{M,r}]}$ is nonquasianalytic, see Lemma \ref{prnonquasi}).

Conversely, if $M$ does not satisfy \hyperlink{mnqr}{$(\on{nq}_r)$},
by  Lemma \ref{prnonquasi} the class $\mathcal{E}_{[P^{M,r}]}$ is quasianalytic and the conclusion follows from  \eqref{spaceinclusion}.
\qed\enddemo

\textbf{Acknowledgements:} The authors wish to thank the anonymous referees for their valuable suggestions
that improved the presentation of this paper. The first two authors are partially supported by the Spanish Ministry of Economy, Industry and Competitiveness under the project MTM2016-77642-C2-1-P. The first author has been partially supported by the University of Valladolid through a Predoctoral Fellowship (2013 call) co-sponsored by the Banco de Santander. The third author is supported by FWF-Project J~3948-N35, as a part of which he has been an external researcher at the Universidad de Valladolid (Spain) for the period October 2016 - December 2018.\par

\bibliographystyle{plain}
\bibliography{Bibliography}

\vskip1cm

\textbf{Affiliation}:\\
J.~Jim\'{e}nez-Garrido, J.~Sanz:\\
Departamento de \'Algebra, An\'alisis Matem\'atico, Geometr{\'\i}a y Topolog{\'\i}a, Universidad de Valladolid\\
Facultad de Ciencias, Paseo de Bel\'en 7, 47011 Valladolid, Spain.\\
Instituto de Investigaci\'on en Matem\'aticas IMUVA\\
E-mails: jjjimenez@am.uva.es (J.~Jim\'{e}nez-Garrido), jsanzg@am.uva.es (J. Sanz).
\\
\vskip.5cm
G.~Schindl:\\
Fakult\"at f\"ur Mathematik, Universit\"at Wien, Oskar-Morgenstern-Platz~1, A-1090 Wien, Austria.\\
E-mail: gerhard.schindl@univie.ac.at.

\end{document}